\documentclass[gdplain]{geradwp}

\PassOptionsToPackage{hyphens}{url}

\usepackage{algorithm2e}
\usepackage{hyperref}

\graphicspath{{Figures/}} 
\hypersetup{colorlinks,%
citecolor={blue}, 
urlcolor={blue},
linkcolor={blue},
breaklinks={true}
}

\usepackage{xspace}
\usepackage{amsmath}
\usepackage{dsfont}
\usepackage{amssymb}
\usepackage{stmaryrd}
\usepackage{mathtools}
\usepackage{subfigure}
\usepackage{soul}
\usepackage{xfrac}
\usepackage[square,numbers]{natbib}

\DeclareMathOperator*{\minimize}{minimize}
\newcommand{\dsum}[2]{\sum\limits_{#1}^{#2}}
\newcommand{\abs}[1]{\left\lvert{#1}\right\rvert}
\newcommand{\norm}[1]{\lVert#1\rVert}
\def\ll{\llbracket}
\def\rr{\rrbracket}
\def\N{\mathbb{N}}
\def\Z{\mathbb{Z}}
\def\R{\mathbb{R}}
\def\1{\mathds{1}}
\newcommand{\e}[1]{e^{#1}}
\renewcommand{\exp}[1]{exp\left(#1\right)}
\newcommand{\floor}[1]{\left\lfloor#1\right\rfloor}
\newcommand{\solver}[1]{\textsf{\textsc{#1}}\xspace}
\newcommand{\algofield}[1]{\textsc{#1}\xspace}
\def\mads{\algofield{Mads}}
\def\nomad{\solver{Nomad}}
\def\ipopt{\solver{Ipopt}}
\def\jump{\solver{JuMP}}
\def\dfo{\algofield{DFBBO}}
\def\so{\algofield{SO}}
\def\P{\mathcal{P}}
\def\Rcal{\mathcal{R}}
\def\S{\mathcal{S}}

\newcommand{\fct}[5]{
    #1:\left\{\begin{array}{ccl}
        #3 & \to     & #5\\
        #2 & \mapsto & #4
    \end{array}\right.}
\newcommand{\fctoneline}[5]{#1:#2\in#3 \mapsto #4\in#5}

\newcommand{\problemoptimcoresup}[3]{
    \begin{array}{cl}
        \underset{
            \begingroup\scriptsize\renewcommand{\arraystretch}{1.25}\begin{array}{c}#2\end{array}\endgroup
        }{#1} & #3
    \end{array}
}
\newcommand{\problemoptimcore}[5]{
    \begin{array}{cl}
        \underset{
            \begingroup\scriptsize\renewcommand{\arraystretch}{1.25}\begin{array}{c}#2\end{array}\endgroup
        }{#1} & \begin{array}{l}#3\end{array}
    \\[#5]
        \mbox{subject to} & \begin{array}[t]{ll}#4\end{array}
    \end{array}
}
\newcommand{\problemoptim}[4]{\problemoptimcore{#1}{#2}{#3}{#4}{1ex}}
\newcommand{\problemoptimcontrol}[4]{\problemoptimcore{#1}{#2}{#3}{#4}{5.5ex}}
\newcommand{\problemoptimoneline}[3]{\problemoptimcoresup{#1}{#2}{#3}}

\newtheorem{theorem}{Theorem}[section]
\newtheorem{remark}[theorem]{Remark}
\newtheorem{definition}[theorem]{Definition}
\newtheorem{proposition}[theorem]{Proposition}
\newtheorem{example}[theorem]{Example}

\GDtitre{A derivative-free approach to optimal control problems with a piecewise constant Mayer cost function}
\GDmois{Novembre}{November}
\GDannee{2021}
\GDnumero{XX}
\GDauteursCourts{C. Audet, P.-Y. Bouchet, L. Bourdin}
\GDauteursCopyright{Bouchet, Audet, Bourdin}
\GDpostpubcitation{Hamel, Benoit, Karine H\'ebert (2021). \og Un exemple de citation \fg, \emph{Journal of Journals}, vol. X issue Y, p. n-m}{https://www.gerad.ca/fr}
\GDrevised{Mai 2021}

\begin{document}

\GDpageCouverture

\begin{GDpagetitre}

\begin{GDauthlist}
\GDauthitem{Pierre-Yves Bouchet \ref{affil:gerad}\GDrefsep\ref{affil:polymtl}}
\GDauthitem{Charles Audet \ref{affil:gerad}\GDrefsep\ref{affil:polymtl}}
\GDauthitem{Loïc Bourdin \ref{affil:unilim}}
\end{GDauthlist}

\begin{GDaffillist}
\GDaffilitem{affil:gerad}{GERAD, Montr\'eal (Qc), Canada, H3T 1J4}
\GDaffilitem{affil:polymtl}{Polytechnique Montr\'eal, Montr\'eal (Qc), Canada, H3T 1J4}
\GDaffilitem{affil:unilim}{XLIM Research Institute, UMR CNRS 7252, University of Limoges, France.}
\end{GDaffillist}

\begin{GDemaillist}
\GDemailitem{pierre-yves.bouchet@polymtl.ca}
\GDemailitem{charles.audet@gerad.ca}
\GDemailitem{loic.bourdin@unilim.fr}
\end{GDemaillist}

\end{GDpagetitre}


\GDabstracts

\begin{GDabstract}{Abstract}
    A piecewise constant Mayer cost function is used to model optimal control problems in which the state space is partitioned into several regions, each having its own Mayer cost value.
    In such a context, the standard numerical methods used in optimal control theory naturally fail, due to the discontinuities and the null gradients associated with the Mayer cost function.
    In this paper an hybrid numerical method, based on both derivative-free optimization and smooth optimization techniques, is proposed to solve this class of problems.
    Numerical simulations are performed on some standard control systems to show the efficiency of the hybrid method, where \nomad and \ipopt are used as, respectively, derivative-free optimization and smooth optimization solvers.

\end{GDabstract}

\begin{GDabstract}{R\'esum\'e}
    Une fonction de coût de Mayer constante par morceaux est requise pour correctement modéliser des problèmes de contrôle optimal dans lesquels l'espace des états est partitonné en différentes régions ayant chacune une valeur de coût de Mayer.
    Dans un tel contexte, les méthodes numériques communes dans la théorie du contrôle optimal échouent naturellement, à cause des discontinuités et des gradients nuls associés à la fonction de coût de Mayer.
    Dans cet article, nous proposons une méthode numérique hybride, exploitant des techniques d'optimisation lisse et d'optimisation sans dérivées, pour traiter cette classe de problèmes.
    Nous présentons des simulations numériques sur plusieurs systèmes de contrôle classiques pour montrer l'efficacité de notre méthode hybride, en exploitant les solveur \ipopt pour l'optimisation lisse et \nomad pour l'optimisation sans dérivées.

\paragraph{Mots cl\'es\,: }
    contrôle optimal,
    fonction de coût de Mayer constante par morceaux,
    optimisation sans dérivées et de boîte noire.

\end{GDabstract}

\begin{GDacknowledgements}
Work of the first author is supported by NSERC Canada Discovery Grant 2020-04448.

\end{GDacknowledgements}


\GDarticlestart
\section{Introduction}
\label{sec:intro}

\subsection{Motivation}
\begin{figure}[h]
    \centering
    \includegraphics[width=0.7\linewidth]{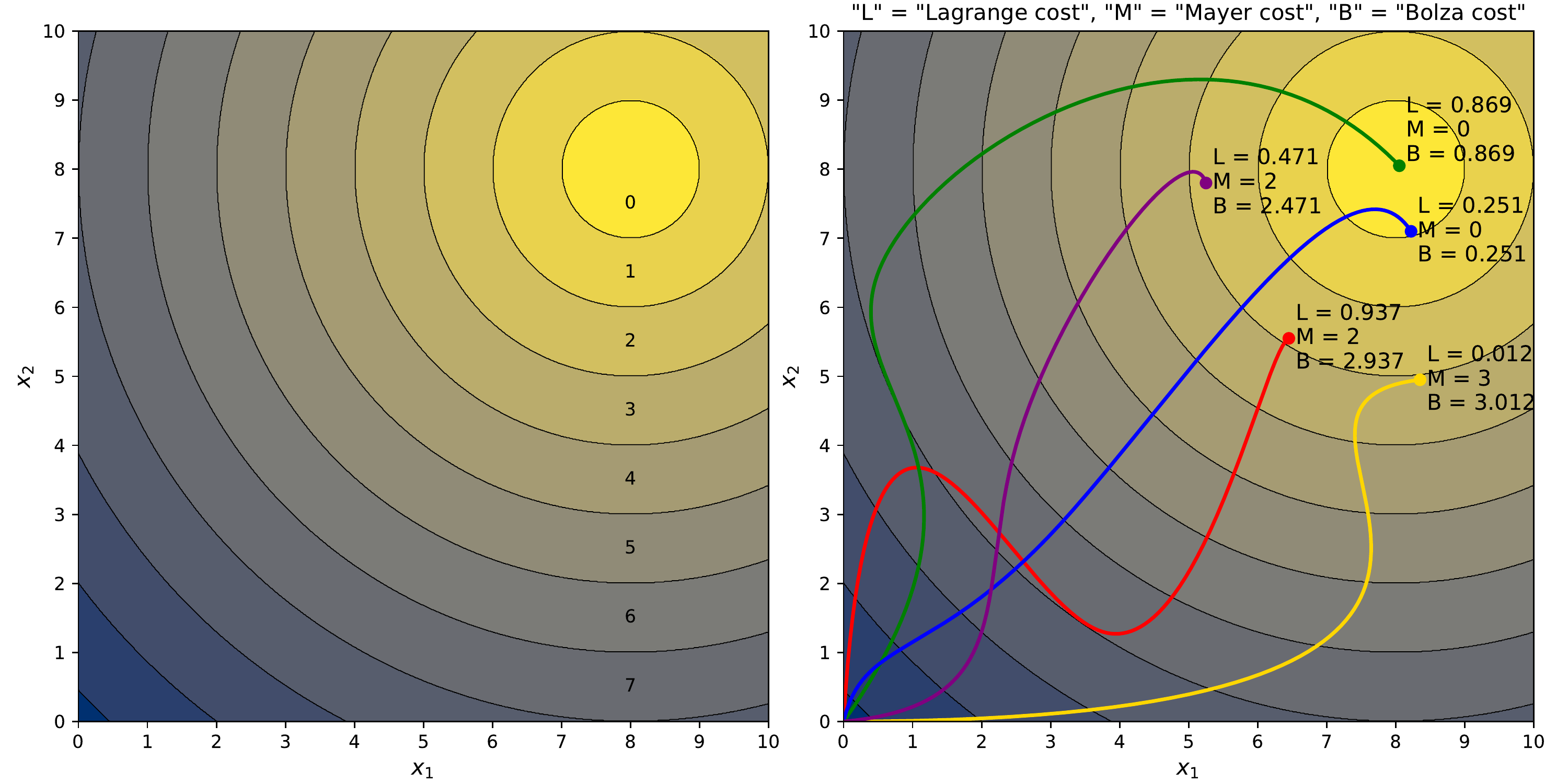}
    \caption{(a) Representation of the Mayer cost function of Problem~\eqref{eq:Pcontrol_continuous_motivation} in the two-dimensional case~$n=2$ with~$x_\mathrm{end} = (8,8)$.
        (b) Artificial trajectories and their associated costs.
        }
    \label{fig:intro_figure}
\end{figure}

The present work was initially motivated by numerically solving Bolza optimal control problems of the form
\begin{equation}
    \tag{\mbox{$\mathrm{P}_\mathrm{Motiv}$}}\label{eq:Pcontrol_continuous_motivation}
    \problemoptimcontrol
        {\minimize}
        {x \in \mathcal{AC} ([0,T],\R^n) \\
         u \in \mathcal{L}^{\infty}([0,T],\R^m)
        }
        {\floor{\norm{x(T)-x_{\mathrm{end}}}}
            +
        \displaystyle \int_0^T \ell(t,x(t),u(t)) \, dt 
        }
        {\dot{x}(t) = f(t,x(t),u(t)), & \mbox{a.e.\ } t \in [0,T], \\[3pt]
        x(0)=x_0, \\[3pt]
        c(t,x(t),u(t)) \leq 0, & \mbox{a.e.\ } t \in [0,T], \\[3pt]
        c_T(x(T)) \leq 0
        }
\end{equation}
where standard notations are used (recalled in Sections~\ref{sec:app_opt_control} and~\ref{sec:numerical_tests})
and in which the dynamics, the constraint functions and the Lagrange cost function are smooth, but where the Mayer cost is discontinuous, being defined as the truncation (involving the floor function~$\floor{\cdot} : \R \to \Z$) of the distance between the final state~$x(T)$ and a reference target~$x_{\mathrm{end}}$. Note that the corresponding Mayer cost function is piecewise constant with integer values. In this paper we call \textit{plateaus} the subsets of~$\R^n$ on which a piecewise constant function from~$\R^n$ to~$\R$ has a given value.
As illustrated in Figure~\ref{fig:intro_figure}(a), in the two-dimensional case $n=2$ the plateaus of the Mayer cost function of Problem~\eqref{eq:Pcontrol_continuous_motivation} are concentric rings centered at~$x_{\mathrm{end}}$.
Such a Mayer cost function is used to model practical problems in which a kind of \textit{discrete distance} to a target is evaluated, for example in precision sports (shooting, darts, archery, precision landing of a parachute, etc.).
Figure~\ref{fig:intro_figure}(b) illustrates several properties of these problems.
For example, two trajectories ending on the same plateau are only discriminated by the Lagrange cost.
Also, a trajectory with a high Lagrange cost ending on a given plateau may have a smaller Bolza cost than another with a low Lagrange cost but ending on a higher-value plateau.
Finally, the two trajectories ending on the lowest-cost plateau show that reaching $x_{\mathrm{end}}$ as close as possible may not be the optimal solution.
These examples show that one cannot remove the truncation in the Mayer cost function without altering (possibly significantly) the optimal solutions.

In optimal control theory, the numerical methods are usually classified as either direct or indirect methods~\cite{RaoNumericalOptControlSurvey,Trelat}.
An indirect method relies on the first-order necessary optimality conditions provided by the \textit{Pontryagin maximum principle} (PMP).
Most of the literature on PMP concerns smooth problems (see, \textit{e.g.},~\cite{BressanPiccoli,BrysonHo,Pontryagin1962}), but extensions exist for problems satisfying (only) Lipschitz-continuity assumptions, by using generalized notions of gradient (see, \textit{e.g.},~\cite{ClarkeBook,FrankowskaNecessary,Vinter2010Book}).
Some papers derive a PMP for problems with discontinuities, such as~\cite{HaberkornTrelat2011} in the context of hybrid dynamics, or~\cite{BayenPfeiffer} in the context of time crisis problems where discontinuities occur as an indicator function in the Lagrange cost function.
However, deriving a relevant version of the PMP in our context is unlikely because optimal trajectories are inclined to end on the boundaries of the plateaus of the Mayer cost function, where gradients and generalized gradients are not defined.
Our numerical simulations provided in Section~\ref{sec:numerical_tests} confirm this phenomenon.
Therefore, in this paper, we rather focus on direct methods which rely on a full discretization of the optimal control problem, via a Runge-Kutta method~\cite{SchwartzPolakRKAnalysis} for example, in order to recover a finite-dimensional optimization problem that can be numerically solved using a relevant optimization algorithm.

However, a full discretization of Problem~\eqref{eq:Pcontrol_continuous_motivation} produces a finite-dimensional optimization problem which may preserve the discontinuity of the Mayer cost function, while most optimization algorithms in the literature are dedicated to smooth problems (see, \textit{e.g.},~\cite{NoWr2006}).
Nevertheless, one can easily design numerical approaches attempting to overcome this obstacle.
A first approach consists in considering a smooth approximation of the discretized problem that can be numerically solved using a \textit{smooth optimization} (\so) solver such as \ipopt~\cite{ipopt}.
However smoothing the piecewise constant Mayer cost function leads to a smooth function with gradients which are null (or almost null) on most of the domain.
Hence they do not provide useful information to the \so algorithm which thus focuses only on the optimization of the Lagrange cost, while the latter may be less rewarding than optimizing the Mayer cost, as illustrated in Figure~\ref{fig:intro_figure}(b).
Our numerical results provided in Section~\ref{sec:numerical_tests} confirm that this first approach is inefficient.
A second approach relies in solving directly the full discretization of Problem~\eqref{eq:Pcontrol_continuous_motivation}, but with a \textit{derivative-free and blackbox optimization} (\dfo) algorithm.
Indeed these algorithms are designed to solve singular finite-dimensional optimization problems.
\dfo relies only on a proximal analysis of candidate solutions close to the current incumbent, and requires neither existence nor values of gradients~\cite{AuHa2017,ConnScheinbergVicente}.
In this paper we focus on the \mads algorithm~\cite{AuDe2006} implemented in the \nomad solver~\cite{Nomad},
but other techniques exist (see, \textit{e.g.},~\cite{AudetConnLeDigabelPeyrega,bagirov2020numerical,CustodioScheinbergVicente,LiuzziLucidiRinaldiVicente}).
\mads can be applied to a full discretization of Problem~\eqref{eq:Pcontrol_continuous_motivation}, even if the discontinuities of the original problem are preserved.
However \mads is usually efficient for problems with a few dozens of variables at most~\cite[Section 1.4]{AuHa2017}, while a discretized optimal control problem usually has a much larger scale.
Our numerical results provided in Section~\ref{sec:numerical_tests} confirm that this second approach is also inefficient.

\subsection{Contributions}

\begin{figure}[h]
    \centering
    \includegraphics[width=0.85\linewidth]{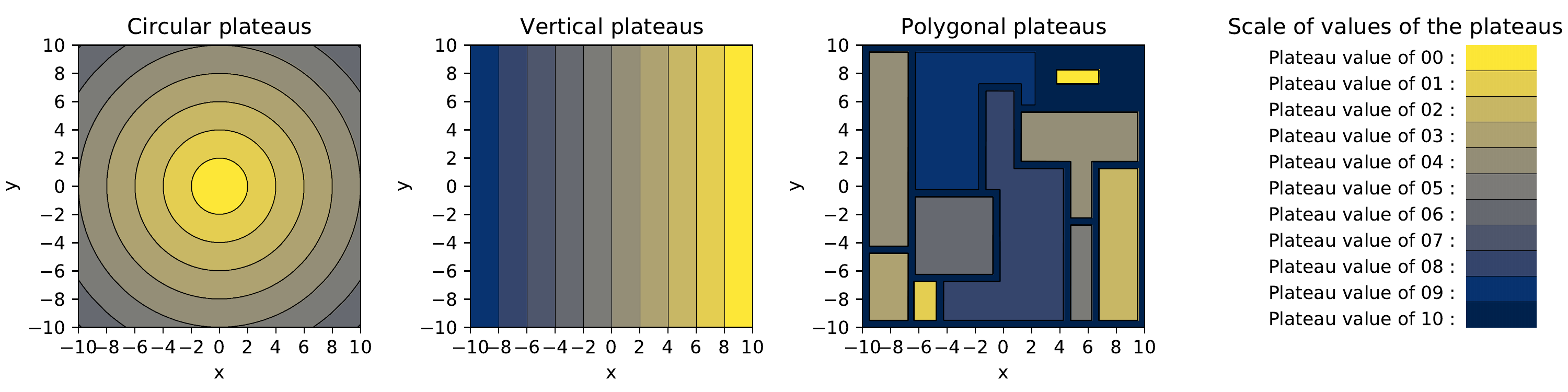}
    \caption{Examples of possible plateaus in the two-dimensional case~$n=2$.
        In all figures of the present paper, the light-colored plateaus represent the lowest-cost plateaus to reach, while the dark ones represent the highest-cost plateaus to avoid.}
    \label{fig:discrete_targets}
\end{figure}

The difficulties mentioned above are actually present in a whole class of optimal control problems.
For example, any piecewise constant Mayer cost function leads to a lack of information provided to \so algorithms, because its gradients are zero (and remain almost-zero after smoothing).
As discussed above, these difficulties cannot be addressed easily with common techniques from the literature.
Therefore, the objective of the present paper is to propose a numerical method to solve efficiently optimal control problems with a piecewise constant Mayer cost function, in view of handling various situations of possible plateaus (as illustrated in Figure~\ref{fig:discrete_targets} in the two-dimensional case~$n=2$).
To illustrate polygonal plateaus as in Figure~\ref{fig:discrete_targets}(c), consider parking a car downtown where we look for a cheap parking lot, while the position of the exact parking spot within this parking lot is less relevant.
In such a problem, the parking lots are the plateaus, each having its own Mayer cost (for example the price to park in the lot), while all parking spots in the same parking lot have the same quality.
However note that reaching a parking spot far from the entrance requires more displacement time.
Hence some local considerations (taken into account by the integral Lagrange cost) may discriminate two parking spots, but this optimization has a noticeably smaller impact than the choice of the parking.

To propose a numerical method to solve efficiently optimal control problems with a piecewise constant Mayer cost function, we exploit the fact that the inner difficulty depends entirely on the final state variable only. 
The main idea is thus to optimize the final state alone via a \dfo algorithm, while a smooth subproblem, which consists in solving the initial problem with the final state being fixed and with a minimal Lagrange cost, is solved via a \so algorithm.
The resulting method is called \textit{hybrid} as it combines tools from \dfo and \so.

The idea to separate variables can be extended to handle cases where the final state is not the only variable affecting a singularity.
Consider for example a variant of Problem~\eqref{eq:Pcontrol_continuous_motivation} where the Mayer cost is defined as~$\floor{\norm{x(T)-x_\mathrm{end}}} + \floor{\norm{x(T/2)-x_\mathrm{mid}}}$, requiring for the trajectory to approach a target~$x_\mathrm{mid}$ at the midtime and a target~$x_\mathrm{end}$ at the final time.
Here the idea is to optimize the states~$x(T/2)$ and $x(T)$ via a \dfo algorithm and solve via a \so algorithm two smooth subproblems recovering a trajectory joining respectively~$x(0)$ to~$x(T/2)$ and~$x(T/2)$ to~$x(T)$, with a minimal Lagrange cost. 

Hence the main generic problem that we address in this work is given by
\begin{equation}
    \tag{\mbox{$\mathrm{P}_\mathrm{Generic}$}}
    \label{eq:Pcontrol_continuous_generic}
    \problemoptimcontrol
        {\minimize}
        {x \in \mathcal{AC} ([0,T],\R^n ) \\
         u \in \mathcal{L}^{\infty}([0,T],\R^m)
        }
        {g(\{x(t)\}_{t\in\Omega})
            +
        \displaystyle \int_0^T \ell(t,x(t),u(t)) \, dt 
        }
        {\dot{x}(t) = f(t,x(t),u(t)), & \mbox{a.e.\ } t \in [0,T], \\[3pt]
        x(0)=x_0, \\[3pt]
        c(t,x(t),u(t)) \leq 0, & \mbox{a.e.\ } t \in [0,T], \\[3pt]
        c_T(x(T)) \leq 0
        }
\end{equation}
where the Mayer cost depends on the state of the control system at one or several times (indexed by the set~$\Omega\subset[0,T]$ with finite cardinality) and the corresponding Mayer cost function is piecewise constant in a generic sense, that is, without being restricted to a discrete distance to some reference targets.
Actually, even if we focus on piecewise constant Mayer cost functions in this paper, our work can be considered for any Mayer cost function presenting singularities such as discontinuities (when involving the floor function~$\floor{\cdot}$ for example).

The hybrid method proposed in this paper relies on the knowledge of the state variables which may affect the singularities. Precisely, once Problem~\eqref{eq:Pcontrol_continuous_generic} is fully discretized and seen as a finite-dimensional optimization problem, we reformulate it as a main problem handling the state variables which may affect the singularities, and involving subproblems dealing with all remaining state-control variables.
The main problem is solved using a \dfo algorithm.
Each subproblem's goal is to recover a partial feasible trajectory-control pair joining with a minimal partial Lagrange cost two consecutive states proposed by the main \dfo algorithm.
Since the variables affecting the singularities are fixed, the subproblems are smooth discrete optimal control problems with fixed terminal states, solved efficiently via a \so algorithm.
Hence the hybrid method applies both \dfo and \so algorithms, 

The goal of this paper is to formally introduce the hybrid method.
Section~\ref{sec:generic_optim_problem} develops a mathematical framework to express precisely the reformulation discussed above, but for a generic singular finite-dimensional optimization problem, allowing to discuss properties of the reformulation in a simple and general context.
The specification to discretized optimal control problems is done in Section~\ref{sec:app_opt_control}, with details on how the hybrid method works in that context.
Section~\ref{sec:numerical_tests} proposes numerical experiments showing that the hybrid method is efficient for optimal control problems with a piecewise constant Mayer cost function.
Finally Section~\ref{sec:conclusion} provides an opening discussion on future perspectives and possible improvements of the hybrid method.

Let us observe that the hybrid method has connections with trajectory optimization~\cite{FischPhD,PetersenCornickBauerRehder} in a context of driving an aircraft.
However, in contrary to the present work, these papers require existence and explicit knowledge of first- and second-order information about the system and the costs.
The hybrid method also has connections with bilevel optimization~\cite{Dempe2020,SinhaMaloDebBilevelReview} in which the common framework defines an upper level trying to minimize the optimal value of a lower level parameterized by the upper level.
However the bilevel paradigm involves a conflictual leader-follower structure, thus it differs from this work which simply partitions the variables.

\section{Reformulation of a generic singular optimization problem}
\label{sec:generic_optim_problem}

This paper proposes a numerical approach dedicated to solve a class of discretized optimal control problems.
These problems can be written as finite-dimensional constrained optimization problems of the form
\begin{equation}
    \tag{\mbox{$\P$}}
    \label{eq:P}
    \problemoptim
        {\minimize}
        {x \in \R^n}
        {\varphi(x)}
        {c(x)\leq 0}
\end{equation}
with $n \geq 1$ real variables~$x \in \R^n$, where $\varphi: \R^n \to \R$ is the \textit{objective function}, and with~$d \geq 1$ inequality constraints expressed as a component-wise vectorial inequality~$c(x) \leq 0$ via the \textit{constraint function} $c: \R^n \to \R^d$.
In some cases, the functions~$\varphi$ or~$c$ may be singular (\textit{e.g.}, discontinuous) with only a known subset of variables affecting the singularities.
Exploiting this knowledge, this section focuses on the \textit{reformulation} of~\eqref{eq:P} given by
\begin{equation}
    \tag{\mbox{$\Rcal$}}
    \label{eq:R}
    \problemoptimoneline
        {\minimize}
        {y\in\R^p}
        {\left(
            \problemoptim{\inf}{z\in\R^q}{\varphi(y,z)}{c(y,z) \leq 0}
        \right)}
\end{equation}
where the $n$ real variables~$x \in \R^n$ are partitioned into two disjoint groups, being expressed as $x=(y,z)$ and ordered so that $y \in \R^p$ stands for the $p \in \ll 0,n \rr$ variables which may affect the singularities of~\eqref{eq:P} (called \textit{singular variables}) and $z \in \R^q$ stands for the remaining~$q = n-p$ variables (called \textit{smooth variables}).
The content in the large parenthesis of~\eqref{eq:R}, depending on $y \in \R^p$, represents the (possibly extended-) real value~$\inf\{\varphi(y,z) : z \in \R^q, \; c(y,z)\leq0\}$.
The format includes situations where~$p=0$ or~$p=n$, but the case of interest is~$p \in \ll 1,n-1 \rr$ with~$n \geq 2$.
In view of analysing Reformulation~\eqref{eq:R}, we introduce for any~$y \in \R^p$ the \textit{subproblem}~\eqref{eq:S} given by
\begin{equation}
    \tag{\mbox{$\S(y)$}}
    \label{eq:S}
    \problemoptim
        {\minimize}
        {z \in \R^q}
        {\varphi(y,z)}
        {c(y,z)\leq0.}
\end{equation}

Reformulation~\eqref{eq:R} is introduced to numerically solve Problem~\eqref{eq:P} when the latter has singularities (\textit{e.g.}, discontinuities) affected only by a known subset of variables. In that context, solving~\eqref{eq:P} numerically may be hard since it is a singular optimization problem in which gradients are not available.
This falls into the \dfo framework in which neither explicit knowledge nor existence of derivatives is required.
However a \dfo algorithm is usually less efficient than a \so algorithm using gradients when those are available~\cite[Section 1.4]{AuHa2017}.
Therefore Reformulation~\eqref{eq:R} is relevant since it separates the variables affecting the singularities of~\eqref{eq:P} from the others.
Hence, the main idea of the hybrid method introduced in this paper is to numerically solve~\eqref{eq:R} instead of~\eqref{eq:P}, by optimizing the variables~$y$ at the singular unconstrained upper level via a \dfo algorithm, while the remaining variables~$z$ are optimized by solving~\eqref{eq:S} via a \so algorithm. 
This approach gives as few variables as possible to a \dfo algorithm and leaves all the others to a smooth subproblem solved with a more efficient \so algorithm.

The section is organized as follows. First Section~\ref{sec:generic_optim_problem/eps-feasibility} 
recalls standard notions, such as the definitions of a \textit{feasible} and a \textit{well-defined} optimization problem.
Section~\ref{sec:generic_optim_problem/reformulation_results}
analyzes relationships existing between~\eqref{eq:P} and~\eqref{eq:R}. Section~\ref{sec:generic_optim_problem/numerical_issues} discusses issues arising in the numerical resolution of~\eqref{eq:R} via the hybrid method.

\subsection{Feasible and well-defined optimization problems}
\label{sec:generic_optim_problem/eps-feasibility}

Consider the optimization problem~\eqref{eq:P}.
The points $x \in \R^n$ such that $c(x) \leq 0$ form the \textit{feasible set}.
A feasible point $x^*$ is said to be an \textit{optimal solution} of~\eqref{eq:P} if it satisfies $\varphi (x^*) \leq \varphi(x)$ for all feasible points $x$.
In that case $\varphi(x^*)$ is called the \textit{optimal value} of~\eqref{eq:P}.

A \textit{constraint violation function} for~\eqref{eq:P} quantifies the feasibility or infeasibility of any point~$x \in \R^n$.
Writing $c(x) = (c_1(x),\dots,c_d(x))$ for all $x \in \R^n$, the constraint violation function $h_c$ considered in this work, as in~\cite{FlLe02a}, is given by
\begin{equation*}\label{eq:violation_function}
    \fct{ h_c }{ x }{ \R^n }{ \dsum{k=1}{d}\max \{ c_k(x),0 \}^2}{ \R }
\end{equation*}
in which the square allows to avoid introducing additional singularities.

\begin{definition}
\label{def:vareps}
    For any $0 \leq \varepsilon \leq +\infty$, the \emph{$\varepsilon$-feasible set} of~\eqref{eq:P} is defined by 
    \begin{equation*}\label{eq:varepsilon-feasible_set}
        \Theta_\varepsilon = \{x \in \R^n : h_c(x) \leq \varepsilon\} 
    \end{equation*}
    and Problem~\eqref{eq:P} is said to be $\varepsilon$\emph{-feasible} (denoted by $\varepsilon$-F) if~$\Theta_\varepsilon \neq \emptyset$, while it is said to be $\varepsilon$\emph{-infeasible} otherwise.
    When~\eqref{eq:P} is $\varepsilon$-F for some $0 \leq \varepsilon \leq +\infty$, it is said to be $\varepsilon$\emph{-feasible-undefined} if $\varphi(\Theta_\varepsilon)$ is not bounded below, and $\varepsilon$\emph{-feasible-defined} ($\varepsilon$-FD) if $\varphi(\Theta_\varepsilon)$ is bounded below, and $\varepsilon$\emph{-feasible-well-defined} ($\varepsilon$-FWD) if $\varphi(\Theta_\varepsilon)$ is bounded below and attains its lower bound.
\end{definition}

In the sequel the $\varepsilon$ will be omitted in the acronyms of Definition~\ref{def:vareps} when it equals~$0$ (\textit{e.g.}, a FWD problem is a $0$-FWD problem).
Note that~\eqref{eq:P} has an optimal solution if and only if~\eqref{eq:P} is FWD.

\begin{definition}
    When it exists, the \emph{minimal infeasibility} of~\eqref{eq:P} is defined as the smallest value $0 \leq \varepsilon \leq +\infty$ for which the set~$\Theta_\varepsilon$ is nonempty.
\end{definition}

Note that \eqref{eq:P} has a minimal infeasibility if and only if the optimization problem
\begin{equation}
    \tag{\mbox{$\mathcal{H}$}}
    \label{eq:H}
    \problemoptimoneline
        {\minimize}
        {x\in\R^n}
        {h_c(x)}
\end{equation}
is FWD.
If so, the optimal value of~\eqref{eq:H} equals the minimal infeasibility of~\eqref{eq:P}.

\begin{remark}
Consider Problem~\eqref{eq:P}.
The set $\Theta_0$ is the feasible set, $\Theta_\infty = \R^n$, and there necessarily exists~$\varepsilon \geq 0$ such that~$\Theta_\varepsilon \neq \emptyset$.
However~\eqref{eq:P} may be infeasible with no minimal infeasibility.
For example, with $\fctoneline{ c }{ x }{ \R }{ \e{x} }{ \R }$ and $\varphi \equiv 0$, there is no smallest value~$0 \leq \varepsilon \leq +\infty$ for which $\Theta_\varepsilon \neq \emptyset$.
\end{remark}

The concepts introduced above for Problem~\eqref{eq:P} are extended to Reformulation~\eqref{eq:R} and to Subproblems~\eqref{eq:S}, replacing $h_c$ by $h_{c(y,\cdot)}$, for all $y \in \R^p$.

\subsection{Relationships between Problem~\eqref{eq:P} and Reformulation~\eqref{eq:R}}
\label{sec:generic_optim_problem/reformulation_results}

The equality
\begin{equation*}
    \left(\problemoptim
        {\inf}
        {x \in \R^n}
        {\varphi(x)}
        {c(x)\leq 0}\right)
    \quad =
    \problemoptimoneline
        {\inf}
        {y\in\R^p}
        {\left(
            \problemoptim{\inf}{z\in\R^q}{\varphi(y,z)}{c(y,z)\leq0}
        \right)}
\end{equation*}
holds true, even when the infima are not attained, or even when they are equal to~$\pm\infty$.
Hereafter we investigate when the infima are attained (and thus are minima) and what relations can be obtained between the minimizers.
The next proposition (whose proof follows directly from Definition \ref{def:vareps}) establishes an equivalence result on the FWD nature of~\eqref{eq:P} and the FWD nature of~\eqref{eq:R}, and, if so, derives some relations between their respective optimal solutions.

\begin{proposition}\label{th:equivalence_problems}
The following statements hold:
\begin{enumerate}
    \item
        If~\eqref{eq:P} is FWD, then, for any optimal solution $x^* = (y^*,z^*)$ of~\eqref{eq:P}, $y^*$ is an optimal solution of~\eqref{eq:R} (which is FWD) and $z^*$ is an optimal solution of~$\S(y^*)$ (which is FWD).
    \item
        If~\eqref{eq:R} is FWD and has an optimal solution $y^*$ for which $\S(y^*)$ is FWD, then, for any optimal solution $z^*$ of $\S(y^*)$, $x^* = (y^*,z^*)$ is an optimal solution of~\eqref{eq:P} (which is FWD).
\end{enumerate}
\end{proposition}

The second item of Proposition~\ref{th:equivalence_problems} guarantees that we can recover an optimal solution of~\eqref{eq:P} from an optimal solution $y^*$ of~\eqref{eq:R}, provided that $\S(y^*)$ is FWD.
The next counterexamples illustrate this requirement, with two different obstructions.
In these examples (and the others of this section), note that the singularity affecting the variable $y$ is not mandatory but allows to remain consistent with our setting.

\begin{example}\label{ex:R_fwd_P_fd}
    Consider $\fctoneline{\varphi}{(y,z)}{\R^2}{\e{\abs{y}+z}}{\R}$ and $c \equiv 0$.
    In that context~\eqref{eq:R} is FWD and any~$y \in \R$ is an optimal solution of~\eqref{eq:R}.
    However, for any $y \in \R$,~\eqref{eq:S} is FD but not FWD.
    Note that~\eqref{eq:P} is FD, but is not FWD.
\end{example}

\begin{example}\label{ex:unusable_solution_of_reformulated_problem}
    Consider
        $\fctoneline{ \varphi }{ (y,z) }{ \R^2 }{ \abs{y}\e{z} + \abs{y-1}z^2 }{ \R }$
    and $c \equiv 0$.
  In that context~\eqref{eq:P} is FWD and its unique optimal solution is $x^* = (0,0)$.
  Furthermore~\eqref{eq:R} is FWD and has two optimal solutions $y^*_1 = 0$ and $y^*_2 = 1$.
    Note that $\S(y^*_1)$ is FWD (with $z^*_1 = 0$ an optimal solution), while $\S(y^*_2)$ is FD but not FWD.
\end{example}

%

Proposition~\ref{th:equivalence_problems} provides only a theoretical result establishing relations between the optimal solution $x^* = (y^*,z^*)$ of~\eqref{eq:P}, the optimal solution $y^*$ of~\eqref{eq:R} and the optimal solution $z^*$ of~$\S(y^*)$.
These relations rely in particular on the FWD nature of $\S(y^*)$.
Note that, in some cases, \eqref{eq:S} may be FWD if and only if~$y=y^*$ (see Example~\ref{ex:not_well_defined_subproblem}).
It is also possible that~\eqref{eq:S} is infeasible for some~$y \in \R^p$ arbitrarily close to~$y^*$ (see Example~\ref{ex:infeasible_subproblem}).

\begin{example}
\label{ex:not_well_defined_subproblem}
    Consider
        $\fctoneline{ \varphi }{ (y,z) }{ \R^2 }{ \abs{y}\e{z} }{ \R }$
    and
        $c \equiv 0$.
   In that context~\eqref{eq:P} has an optimal solution $x^* = (y^*,z^*) = (0,0)$.
    However, when $y \neq y^*$, $\S(y)$ is FD but is not FWD.
\end{example}


\begin{example}\label{ex:infeasible_subproblem}
    Consider
        $\fctoneline{ \varphi }{ (y,z) }{ \R^2 }{ \abs{y} + z^2 }{ \R }$
    and
        $\fctoneline{ c }{ (y,z) }{ \R^2 }{ (1-y-z,1-y+z) }{ \R^2 }$.
   In that context~\eqref{eq:P} has an optimal solution $x^* = (y^*,z^*) = (1,0)$, while~\eqref{eq:S} is infeasible for any~$y < 1$.
\end{example}

Examples~\ref{ex:not_well_defined_subproblem} and~\ref{ex:infeasible_subproblem} show that, when looking for a (numerical) optimal solution of~\eqref{eq:R}, we may be confronted to subproblems~\eqref{eq:S} that are FD but not FWD, or even infeasible.
These situations have to be addressed carefully because any attempt to evaluate the objective value of~\eqref{eq:R} at some~$y \in \R^p$ requires, roughly speaking, a (numerical) optimal value of~\eqref{eq:S}.
A FD subproblem~\eqref{eq:S} is not hard to handle, as we may obtain a feasible point~$z \in \R^q$ with objective value arbitrarily close to the infimum, but an infeasible subproblem~\eqref{eq:S} is more problematic.
The next section discusses a dedicated numerical process to face such an infeasibility.
It relies on the next proposition, whose proof is simple and left to the reader.

\begin{proposition}\label{th:subproblem_defined}
If~$h_{c(y,\cdot)}$ admits a minimum for all $y \in \R^p$, then $\S(y)$ has a minimal infeasibility, denoted by~$\mu_y$, for all $y \in \R^p$.
If moreover~\eqref{eq:P} is~$\infty$-FD, then~$\S(y)$ is $\mu_y$-FD for all~$y \in \R^p$.
\end{proposition}


\subsection{Issues concerning the numerical resolution of Reformulation~\eqref{eq:R}}
\label{sec:generic_optim_problem/numerical_issues}

In the hybrid method, the objective function of~\eqref{eq:R} is treated as a \textit{blackbox}, that is, as a function with intractable expression which requires here, roughly speaking, to solve a subproblem~\eqref{eq:S} depending on the input $y \in \R^p$.
This blackbox has to be handled carefully since two main issues (discussed and illustrated in Section~\ref{sec:generic_optim_problem/reformulation_results}) need to be addressed.
First~\eqref{eq:S} may be infeasible and, second, denoting by $\mu_y$ its minimal infeasibility (when it exists), \eqref{eq:S} may be~$\mu_y$-FD but not~$\mu_y$-FWD.
Assuming that the conditions expressed in Proposition~\ref{th:subproblem_defined} are satisfied, the following two-phases numerical process is proposed to handle these drawbacks:

\begin{itemize}
    \item[\rm{(i)}]
        Solve with a \so algorithm the smooth unconstrained FWD problem
        \begin{equation}
            \tag{\mbox{$\mathcal{H}_y$}}
            \label{eq:Hy}
            \problemoptimoneline
                {\minimize}
                {z \in\R^q}
                {h_{c(y,\cdot)}(z).}
        \end{equation}
        Denote by $\tilde{z}_y$ a (numerical) optimal solution of~\eqref{eq:Hy} and get $\tilde{\mu}_y = h_{c(y,\cdot)}(\tilde{z}_y)$ as (numerical) approximation of the minimal infeasibility of~\eqref{eq:S}.
    \item[\rm{(ii)}]
        Solve with a \so algorithm the smooth constrained FD problem
        \begin{equation}
            \tag{\mbox{$\tilde{\S}_y$}}
            \label{eq:Stilde}
            \problemoptim
                {\minimize}
                {z \in \R^q}
                {\varphi(y,z)}
                {h_{c(y,\cdot)} (z) \leq \tilde{\mu}_y}
        \end{equation}
        from the starting point $\tilde{z}_y$.
        Denote by $\tilde{z}_y^*$ the (numerical) optimal solution of~\eqref{eq:Stilde} and get~$\varphi(y,\tilde{z}_y^*)$ as (numerical) approximation of the optimal value of~\eqref{eq:Stilde}.
        In practice it is not restrictive that~\eqref{eq:Stilde} may be not FWD, since it is FD and solved up to some numerical thresholds.
        Return $\tilde{\mu}_y$ and $\varphi(y,\tilde{z}_y^*)$.
\end{itemize}

The goal of the above two-phases process is the following.
If~\eqref{eq:S} is feasible for a given~$y\in\R^p$, then the first step leads to a numerically feasible solution of~\eqref{eq:S}, the second step numerically solves~\eqref{eq:S} and the process returns a numerically optimal value of~\eqref{eq:S}.
Otherwise, if~\eqref{eq:S} is infeasible, then an arbitrary choice has to be made to assign a value to the blackbox at~$y$.
We could set the value $+\infty$, but such a choice would not provide any information to pursue efficiently the numerical resolution of~\eqref{eq:R}.
Instead, via the above process, we recover an infeasibility measure, as well as a corresponding finite objective value.
Recall that \dfo algorithms can exploit these two (possibly conflictual) metrics.
For example the \mads algorithm~\cite{AuDe2006} implemented in \nomad~\cite{Nomad} handles them via its \textit{progressive barrier} formulation~\cite{ProgressiveBarrier}.
It constructs an approximated Pareto front of undominated infeasible solutions on a diagram showing objective value versus infeasibility.
The front is truncated to eliminate solutions with an infeasibility above a threshold decreasing with the iterations.
Each iteration of \mads starts from the solution with the lowest objective function value among the truncated front, so that the algorithm eventually finds a feasible solution with a low objective value.
Example~\ref{ex:convention} illustrates the relevancy of this approach.

\begin{example}\label{ex:convention}
    Consider Example~\ref{ex:infeasible_subproblem} and an attempt to numerically solve~\eqref{eq:R} from the initial guess~$y=0$.
    Subproblem~\eqref{eq:S} is infeasible for any~$y < 1$ and the above numerical two-phases process returns~$\tilde{\mu}_y = 2(1-y)^2$ and~$\varphi(y,\tilde{z}^*_y) = y$.
    Hence any value~$0 < y \leq 1$ has an higher objective function value than $y=0$ but a lower infeasibility.
    An algorithm solving~\eqref{eq:R} can be driven to~$y=1$, leading to a feasible subproblem~\eqref{eq:S}, via an iterative reduction of the maximal admissible infeasibility.
    Then it eventually reaches $x = (y,\tilde{z}^*_{y}) = (1,0)$ which is the optimal solution of~\eqref{eq:P}.
\end{example}


\section{Numerical setup to solve Problem~\eqref{eq:Pcontrol_continuous_generic}}
\label{sec:app_opt_control}

A full discretization of Problem~\eqref{eq:Pcontrol_continuous_generic}, using a Runge-Kutta method for example, can be written as
\begin{equation}
    \tag{\mbox{$\mathrm{P}$}}
    \label{eq:Pcontrol}
    \problemoptimcontrol{\minimize}
        {\{x_k\}_{k=1}^{N} \subset \R^n,\\
         \{u_k\}_{k=0}^{N-1} \subset \R^m
        }
        {g\left(\{x_k\}_{k\in\Omega}\right) + \dsum{k=0}{N-1}\ell_k\left(x_k,u_k\right)}
        {x_{k+1}        =  f_k(x_k,u_k), & k \in \ll0,N-1\rr,\\
         c_k(x_k,u_k) \leq 0,            & k \in \ll0,N-1\rr, \\
         c_N(x_N)     \leq 0
        }
\end{equation}
where
    $n \geq 1$,
    $m \geq 1$ and $N \geq 1$ are fixed positive integers,
    $\{x_k\}_{k=1}^{N} \subset \R^n$ are the \textit{state variables} and
    $\{u_k\}_{k=0}^{N-1} \subset \R^m$ are the \textit{control variables},
    $x_0 \in \R^n$ is the fixed initial state,
    $\Omega \subseteq \ll1,N\rr$ is the set (with $\omega \geq 0$ elements) of the indices $k$ such that the piecewise constant \textit{Mayer cost function} $g : (\R^n)^\omega \to \R$ actually depends on~$x_k$,
    each~$\ell_k : \R^n \to \R$ is the \textit{Lagrange cost function at time~$k$},
    each~$f_k : \R^n \times \R^m \to \R^n$ is the \textit{dynamics at time~$k$},
    each~$c_k : \R^n \times \R^m \to \R^{d_k}$ (with~$d_k \geq 1$) is the \textit{mixed state-control constraints function at time~$k$} and
    $c_N : \R^n \to \R^{d_N}$ (with~$d_N \geq 1$) stands for the \textit{final state constraints function}.
    
Note that Problem~\eqref{eq:Pcontrol} is a specific instance of Problem~\eqref{eq:P}.
Therefore, following the approach developed in Section~\ref{sec:generic_optim_problem}, Section~\ref{sec:app_opt_control/reformulated_problem} proposes several possible reformulations of~\eqref{eq:Pcontrol} to handle its singularities.
Then Section~\ref{sec:app_opt_control/numerical_methods} provides a numerical setup to introduce the \dfo-\so hybrid method based on these reformulations, and to compare it with two non-hybrid approaches (one based on a \so algorithm exclusively, and one based on a \dfo algorithm exclusively).
Some numerical experiments, showing the efficiency of the hybrid method for optimal control problems with a piecewise constant Mayer cost function, are provided in the next Section~\ref{sec:numerical_tests}.


%
%

\subsection{Several reformulations of Problem~\eqref{eq:Pcontrol}}
\label{sec:app_opt_control/reformulated_problem}

Take $(i,j) \in \ll0,N\rr^2$ with~$j > i$ and $(x_i,x_j) \in (\R^n)^2$.
The search of a trajectory linking~$x_i$ at time $i$ to~$x_j$ at time~$j$, with minimal partial Lagrange cost and satisfying the state-control constraints of Problem~\eqref{eq:Pcontrol}, leads to the subproblem given by
\begin{equation}
    \tag{\mbox{$\mathrm{S}(i,x_i,j,x_j)$}}
    \label{eq:Scontrol}
    \problemoptimcontrol{\minimize}
        {\{x_k\}_{k=i+1}^{j-1} \subset \R^n,\\
         \{u_k\}_{k=i\phantom{+1}}^{j-1} \subset \R^m}
        {\dsum{k=i}{j-1} \ell_k\left(x_k,u_k\right)}
        {x_{k+1} = f_k(x_k,u_k), & k \in \ll i,j-1\rr, \\
         c_k(x_k,u_k) \leq 0,    & k \in \ll i,j-1\rr.
        }
\end{equation}
This subproblem is not always FWD in general.
Hereafter the (possibly extended-real) lower bound of its objective function over its feasible set is denoted by $\underline{S}(i,x_i,j,x_j)$.

Now, following the approach developed in Section~\ref{sec:generic_optim_problem} and using the notation introduced above, Problem~\eqref{eq:Pcontrol} has the reformulation 
\begin{equation}
    \tag{\mbox{$\mathrm{R}^{\ll1,N\rr}$}}
    \label{eq:Rcontrol_allvar}
    \problemoptim{\minimize}
        {\{x_k\}_{k=1}^{N} \subset \R^n}
        {g\left(\{x_k\}_{k\in\Omega}\right) + \dsum{k=0}{N-1} \underline{S}(k, x_k, k+1, x_{k+1})}
        {c_N(x_N) \leq 0}
\end{equation}
where the variables are the entire trajectory, instead of the entire trajectory-control pair in Problem~\eqref{eq:Pcontrol}.
For any given trajectory $\{x_k\}_{k=1}^{N}$, the search of a control $\{u_k\}_{k=0}^{N-1}$ driving $x_0$ along this trajectory is implicit.
The search of each control~$u_k$ linking~$x_k$ to~$x_{k+1}$, while minimizing $\ell_k(x_k,u_k)$ and satisfying the constraint~$c_k(x_k,u_k) \leq 0$, is done in the smooth Subproblem~$(\mathrm{S}(k,x_k,k+1,x_{k+1}))$.
However there is still an important number of variables (precisely~$nN$) left to Reformulation~\eqref{eq:Rcontrol_allvar}.

It is possible to reduce this number of variables. The variant
\begin{equation}
    \tag{\mbox{$\mathrm{R}^{\{N\}}$}}
    \label{eq:Rcontrol_xm}
    \problemoptim{\minimize}
        {x_N \in \R^n}
        {g\left(x_N\right) + \underline{S}\left(0, x_0, N, x_N\right)}
        {c_N(x_N) \leq 0}
\end{equation}
has only the final state $x_N$ as a variable.
For any given $x_N \in \R^n$, the search of a control linking $x_0$ to~$x_N$, with minimal Lagrange cost and satisfying the state-control constraints of Problem~\eqref{eq:Pcontrol}, is left to the unique subproblem~$(\mathrm{S}(0, x_0, N, x_N))$.

When $N$ is even and $\Omega \subseteq \{N/2,N\}$, another variant in which the variables are the final state $x_N$ and the midterm state $x_{N/2}$ is possible.
In that case there are two subproblems. The first one links $x_0$ to $x_{N/2}$ and the second one links $x_{N/2}$ to $x_N$.
Denoting by $\lambda_0=0$, $\lambda_1=N/2$,~$\lambda_2=N$, the corresponding reformulation is given by
\begin{equation}
    \tag{\mbox{$\mathrm{R}^{\left\{N/2,N\right\}}$}}
    \label{eq:Rcontrol_xmsur2xm}
    \problemoptim{\minimize}
        {\{x_{\lambda_1},x_{\lambda_2}\} \subset \R^n}
        {g\left(\{x_k\}_{k\in\Omega}\right)
            +
            \dsum{i=0}{1} \underline{S}\left(\lambda_i, x_{\lambda_i}, \lambda_{i+1}, x_{\lambda_{i+1}}\right)
        }
        {c_N(x_{\lambda_2}) \leq 0.}
\end{equation}

Finally, from a general point of view, it is possible to choose \textit{a priori} which states remain variables of the reformulated problem. To this aim, take~$r \in \ll1,N\rr$ and~$\Lambda = \{\lambda_1,\dots,\lambda_r\} \subseteq \ll1,N\rr$ the set of all indices of the states given as variables to the reformulated problem, with increasing values and satisfying $\lambda_r=N$ and $\Omega \subseteq \Lambda$.
Taking $\lambda_r=N$ is mandatory because the constraint $c_N(x_N) \leq 0$ cannot appear in the last subproblem, and taking~$\Omega \subseteq \Lambda$ is required because the Mayer cost cannot appear in any subproblem.
Define also for convenience~$\lambda_0=0$. Finally the general reformulation of Problem~\eqref{eq:Pcontrol} is given by
\begin{equation}
    \tag{\mbox{$\mathrm{R}^\Lambda$}}
    \label{eq:Rcontrol_generic}
    \problemoptim{\minimize}
        {\{x_{\lambda_i}\}_{i=1}^{r} \subset \R^n}
        {g\left(\{x_k\}_{k\in\Omega}\right)
            +
            \dsum{i=0}{r-1} \underline{S}\left(\lambda_i, x_{\lambda_i}, \lambda_{i+1}, x_{\lambda_{i+1}}\right)
        }
        {c_N(x_N) \leq 0}
\end{equation}
with~$nr$ variables and admitting $r$ subproblems, each trying to find a control linking~$x_{\lambda_i}$ to~$x_{\lambda_{i+1}}$ with a minimal partial Lagrange cost while satisfying the state-control constraints of Problem~\eqref{eq:Pcontrol}. Note that Reformulations~\eqref{eq:Rcontrol_allvar}, \eqref{eq:Rcontrol_xm} and \eqref{eq:Rcontrol_xmsur2xm} are particular instances of Reformulation~\eqref{eq:Rcontrol_generic}.

The choices of $r \in \ll1,N\rr$ and~$\Lambda = \{\lambda_1,\dots,\lambda_r\} \subseteq \ll1,N\rr$ follow the knowledge of which states~$x_k$ may affect the singularities of Problem~\eqref{eq:Pcontrol}.
Furthermore note that the smaller $r$ is, then the smaller the number of variables left to Reformulation~\eqref{eq:Rcontrol_generic} solved with a \dfo algorithm is.
In that context note that \eqref{eq:Pcontrol} and \eqref{eq:Rcontrol_generic} are specific instances of \eqref{eq:P} and \eqref{eq:R} from Section~\ref{sec:generic_optim_problem}.
Indeed the $p$ singular variables are the state variables indexed by $\Lambda$ and, for any input $y = (x_{\lambda_i})_{i=1}^p$, Reformulation~\eqref{eq:R} has a (unique) Subproblem~\eqref{eq:S} obtained by aggregating all independent subproblems~$(\mathrm{S}(\lambda_i, x_{\lambda_i}, \lambda_{i+1}, x_{\lambda_{i+1}}))$ together.

\subsection{Numerical setup}
\label{sec:app_opt_control/numerical_methods}

In the next Section~\ref{sec:numerical_tests}, numerical simulations will be performed on examples of Problem~\eqref{eq:Pcontrol} in which the Mayer cost function involves the floor function~$\floor{\cdot} $, causing discontinuities.
This section describes three numerical methods that will be implemented and compared.
The first one is a \so method, the second one is a \dfo method and the third one is the \dfo-\so hybrid method based on the reformulations introduced in Section~\ref{sec:app_opt_control/reformulated_problem}.
Before coming to these methods, we introduce first a simple initialization used for each method.

\newcommand{\methodnocontrol}{$\mathcal{T}_{\mathrm{uncontrolled}}$\xspace}
\paragraph{\textbf{Uncontrolled trajectory \methodnocontrol}}
This trajectory results from the induction~$x_{k+1} = f_k(x_k,u_k)$ for all $k \in \ll0,N-1\rr$, with each control~$u_k$ set to zero.
It is used as initialization of the three numerical methods introduced below.

\newcommand{\methodso}{$\mathcal{M}_{\mathrm{SO}}$\xspace}
\paragraph{\textbf{The \methodso method}}
The purpose of the \methodso method is to solve Problem~\eqref{eq:Pcontrol} by using a \so algorithm, via the JuMP modeling language~\cite{JuMP} and the \ipopt solver, despite the discontinuities due to the floor function~$\floor{\cdot}$.
Unfortunately, \ipopt cannot deal with discontinuous problems.
To overcome this obstacle, we use a smooth approximation of the floor function (see Remark~\ref{rq:smoothing_floor} below) which leads to a smoothed version of~\eqref{eq:Pcontrol} accepted by the code.
However the \methodso method remains likely inefficient.
Indeed, even after a smoothing, the plateaus resulting from the floor function remain and lead to partial derivatives close to zero.
Hence, the \methodso method may fail to detect the significant gains on the objective value it would obtain by making some variables cross the floor function plateaus.
The presumed ineffectiveness of the \methodso method is confirmed by the numerical results in Section~\ref{sec:numerical_tests}.

\newcommand{\methoddfbbo}{$\mathcal{M}_{\mathrm{DFBBO}}$\xspace}
\paragraph{\textbf{The \methoddfbbo method}}
The \methoddfbbo method aims to solve Problem~\eqref{eq:Pcontrol} thanks to the \dfo algorithm \mads via the JuMP interface of the \nomad solver~\cite{NomadJl}.
However \nomad is usually restricted to small-scale problems~\cite[Section~1.4]{AuHa2017} while~\eqref{eq:Pcontrol} usually has a large scale.
Fortunately a dimension reduction is possible, as we can replace all state variables~$x_k$ by the trajectory generated by the control variables~$u_k$ and the equality constraints~$x_{k+1} = f_k(x_k,u_k)$.
However, since each control $u_k$ remains a variable, the scale of the reduced problem remains too large for \nomad, in which a critical value of $50$ is indicated as the maximal number of variables that can be easily handled~\cite{Nomad}.
Hence, the \methoddfbbo method is likely inefficient, which is confirmed by the numerical results obtained in Section~\ref{sec:numerical_tests}.
Note that, for fair comparison with the \methodso method, the \methoddfbbo method is executed in Section~\ref{sec:numerical_tests} on the smoothed version of Problem~\eqref{eq:Pcontrol}, in addition to the above dimension reduction.

\newcommand{\methodhybridgeneric}{$\mathcal{M}^{\Lambda }_{\mathrm{hybrid}}$\xspace}
\newcommand{\methodhybridfinal}{$\mathcal{M}^{\{N\} }_{\mathrm{hybrid}}$\xspace}
\newcommand{\methodhybridfinalmidterm}{$\mathcal{M}^{\{N/2,N\}}_{\mathrm{hybrid}}$\xspace}
\paragraph{\textbf{The \methodhybridgeneric method}}
The \methodhybridgeneric method is the main contribution of the present paper.
Instead of numerically solving Problem~\eqref{eq:Pcontrol}, the \methodhybridgeneric method focuses on Reformulation~\eqref{eq:Rcontrol_generic} introduced in Section~\ref{sec:app_opt_control/reformulated_problem} and numerically solves it via \nomad, while its objective function is evaluated by solving the smooth subproblems~($\mathrm{S}(\lambda_i, x_{\lambda_i}, \lambda_{i+1}, x_{\lambda_{i+1}})$) via \ipopt.
This evaluation relies on the numerical two-phases process presented in Section~\ref{sec:generic_optim_problem/numerical_issues}, where the infeasibility and the corresponding finite objective value are returned to \nomad to be used in its implementation of the progressive barrier variant of the \mads algorithm~\cite{ProgressiveBarrier}.
Unlike the \methodso and \methoddfbbo methods, the \methodhybridgeneric applies each solver to a problem it is designed for.
Indeed \nomad solves a small-scale blackbox optimization problem with complex internal structure, while \ipopt solves smooth optimization problems.
Note that, for fair comparison with the two methods above, the \methodhybridgeneric method is executed in Section~\ref{sec:numerical_tests} on the reformulation of a smoothed version of Problem~\eqref{eq:Pcontrol}, and not on the reformulation of~\eqref{eq:Pcontrol}.
Also, as in Section~\ref{sec:numerical_tests} the examples have singularities affected by the state variables $x_N$ and $x_{N/2}$ only, the two variants \methodhybridfinal and \methodhybridfinalmidterm of the \methodhybridgeneric method will be tested.

\begin{remark}\label{rq:smoothing_floor}
As \ipopt cannot handle discontinuities from the floor function~$\floor{\cdot}$, a smooth approximation is introduced as follows.
First, note that the Heaviside function $\1_{\R_+}(x) = 0$ if~$x < 0$ and $\1_{\R_+}(x) = 1$ if~$x \geq 0$, can be approximated by~$\1_{\R^+}(x) \approx \exp{-\e{-x}}$ for any $x \in \R$.
Strengthening the approximation by a factor~$\tau \in \R_+^*$ and shifting the discontinuity from~$0$ to any value $q \in \R$, we obtain
\begin{equation*}
    \forall (\tau,q) \in \R_+^* \times \R, \quad \forall x \in \R, \quad
    \1_{[q,+\infty[}(x) \approx \exp{-\e{-\tau(x-q)}}.
\end{equation*}

Hence a smooth approximation of~$\floor{\cdot}$ can be derived on any interval $[a,b] \subset \R$, with~$(a,b) \in \Z^2$, $a < b$, and with any strengthening factor $\tau \in \R_+^*$, as follows:
\begin{equation}
    \tag{\mbox{$\floor{\cdot}_{\tau;a.b}$}}
    \label{eq:floor_approx}
    \fct{ \floor{\cdot}_{\tau;a,b} }{ x }{ [a,b] }{ a+\dsum{q=a+1}{b-1}\exp{-\e{-\tau(x-q)}}.}{ [a,b] }
\end{equation}

Note that~\ref{eq:floor_approx} is of class ${\cal C}^\infty$ for any admissible couple $(\tau;a,b)$. Its derivative is close to $0$ at any $x$ far from $\Z$.
Also note that~\ref{eq:floor_approx} is not constant on the intervals~$[q,q+1)$ with $q \in \Z$, and it does not reach integer values when $x\in\Z$.
Figure~\ref{fig:approx_cinfty} shows $\floor{\cdot}$, as well as some approximations~\ref{eq:floor_approx} for some values of $\tau$, $a$ and $b$.
The numerical tests presented in Section~\ref{sec:numerical_tests} are performed with $\tau = 75$.
\begin{figure}[h]
    \centering
    \includegraphics[width=0.75\linewidth]{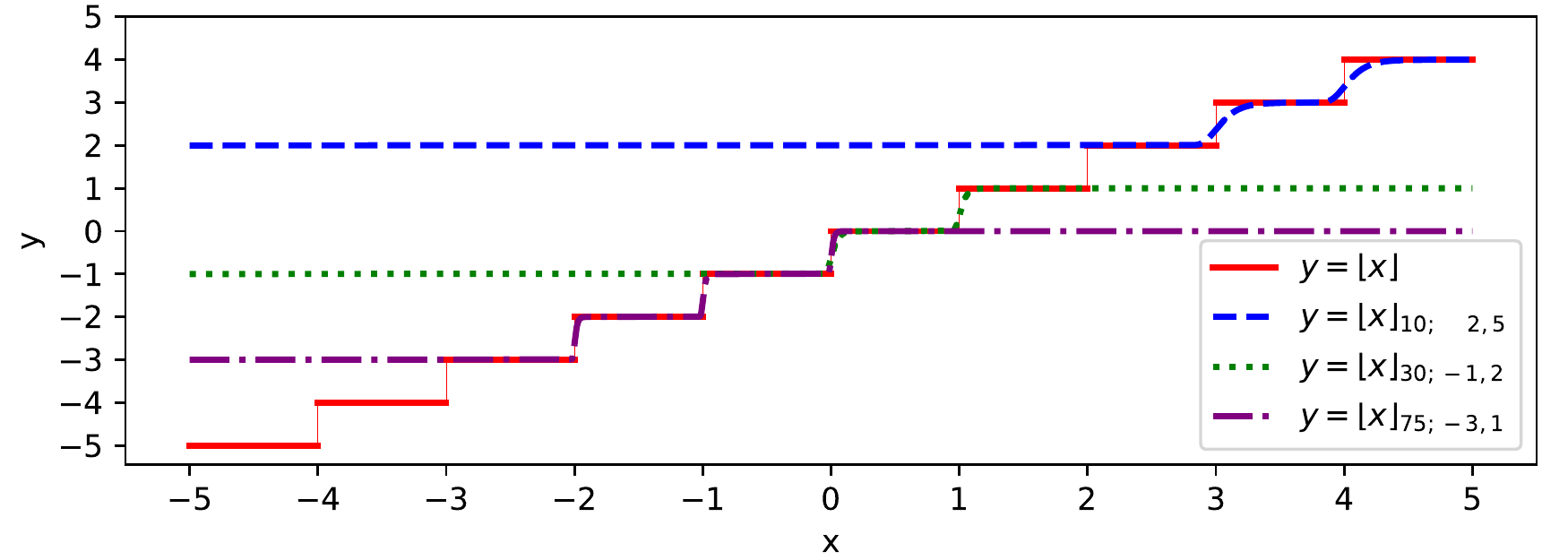}
    \caption{Floor function $\floor{\cdot}$ and smooth approximations~\ref{eq:floor_approx}.}
    \label{fig:approx_cinfty}
\end{figure}

\end{remark}

\section{Computational study}
\label{sec:numerical_tests}

The examples studied in this section are of the form of Problem~\eqref{eq:Pcontrol_continuous_generic} where the Mayer cost function involves the floor function~$\floor{\cdot}$, causing discontinuities.
We discretize the timelapse~$[0,T]$ (with~$T>0$) with a constant timestep~$T/N$ (with $N \in \N^*$), the control system via the standard fourth-order Runge-Kutta method~\cite{SchwartzPolakRKAnalysis} and the integral Lagrange cost via the left rectangle rule.
This discretization procedure leads to Problem~\eqref{eq:Pcontrol} which preserves the singularities due to the floor function.
We implement and compare the numerical methods introduced in Section~\ref{sec:app_opt_control/numerical_methods} in the \href{https://julialang.org/}{\textit{Julia} language} using the modeling tool \jump~\cite{JuMP}.
Both \ipopt and \nomad~\cite{NomadJl} are interfaced with \jump. As explained in Section~\ref{sec:app_opt_control/numerical_methods}, for fair comparisons between the numerical methods, Problem~\eqref{eq:Pcontrol} is preconditioned using the smoothing procedure described in Remark~\ref{rq:smoothing_floor}.

This section is organized as follows.
Section~\ref{sec:numerical_tests/harmonic_oscillator} deals with a controlled harmonic oscillator, where the Mayer cost involves the final kinetic energy minus a truncation of the final potential energy.
Section~\ref{sec:numerical_tests/zermelo_1} discusses a Zermelo-type problem representing a boat on a river with a stream, trying to reach a relevant parking in an harbour.
Section~\ref{sec:numerical_tests/zermelo_2} addresses a similar Zermelo-type problem but with two targets to reach, one at the midtime~$T/2$ and one at the final time~$T$.
Finally Section~\ref{sec:numerical_tests/lotka-volterra} deals with a truncated version of the Lotka-Volterra stabilization problem and shows that the best choice of $\Lambda$ in the \methodhybridgeneric method may be nontrivial.

\paragraph{Notation}
Hereafter $\mathcal{AC}([0,T],\R^n)$ denotes the space of absolutely continuous functions defined on $[0,T]$ with values in $\R^n$ and $\mathcal{L}^\infty([0,T],\R^m)$ denotes the usual Lebesgue space of essentially bounded control functions defined almost everywhere on~$[0,T]$ with values in $\R^m$.
The acronym \textit{a.e.}\ stands for \textit{for almost every}.
In view of using different weighted Euclidean norms in Sections~\ref{sec:numerical_tests/zermelo_2} and~\ref{sec:numerical_tests/lotka-volterra}, the notation $\norm{X}_\alpha$ stands for the Euclidean norm of the vector $(\alpha_1 X_1, \dots, \alpha_n X_n)$ for all~$X$, $\alpha \in \R^n$.

\subsection{Truncated variant of the harmonic oscillator problem}
\label{sec:numerical_tests/harmonic_oscillator}
Consider
\begin{equation}
    \tag{\mbox{$\mathcal{HO}$}}
    \label{pb:harmonic}
    \problemoptimcontrol
        {\minimize}
        {(x,v) \in \mathcal{AC} ([0,T],\R^2 ) \\
         u \in \mathcal{L}^{\infty}([0,T],\R)
        }
        { \displaystyle \frac{M}{2}v(T)^2 - \frac{K}{2} \lfloor x(T) \rfloor^2
            +
        \frac{1}{T} \int_0^T u(t) \, dt 
        }
        {\dot{x}(t) = v(t), & \mbox{a.e.\ } t \in [0,T], \\[3pt]
        \dot{v}(t) = \frac{1}{M} ( - Kx(t)+u(t) ), & \mbox{a.e.\ } t \in [0,T], \\[3pt]
        x(0)=1, \\[3pt]
        v(0)=0, \\[3pt]
        u(t) \in [0,1], & \mbox{a.e.\ } t \in [0,T]
        }
\end{equation}
where the control system mimics a spring-mass system $M\ddot{x}(t) + Kx(t) = u(t)$,~$x(t)$ and~$v(t)$ denote the position and velocity of a point of mass $M>0$ attached to a spring of constant~$K>0$, and where~$u(t) \in [0,1]$ denotes a constrained external force.
Here the objective is to minimize the (smooth) final kinetic energy~$\frac{M}{2}v(T)^2$, plus the average control value~$\frac{1}{T} \int_0^T u(t) \, dt$, minus the (singular) final truncated potential energy~$\frac{K}{2} \lfloor x(T) \rfloor^2$. In this example, note that the Mayer cost function is not piecewise constant, but involves the floor function~$\floor{\cdot}$, causing discontinuities.

\begin{figure}[h]
    \centering
    \includegraphics[width=\linewidth]{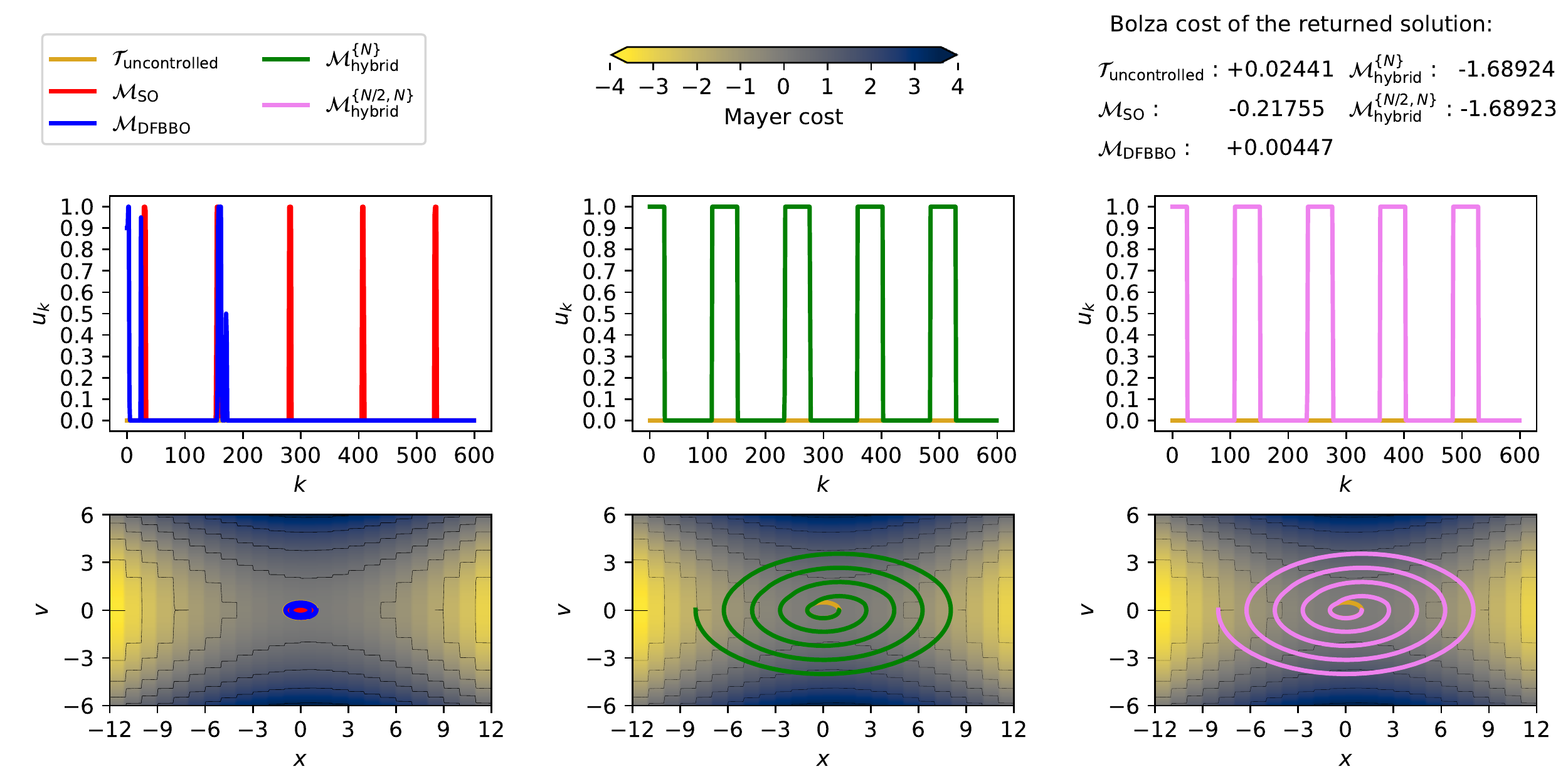}
    \caption{Numerical results obtained on Problem~\eqref{pb:harmonic}.
        The~\methodso and~\methoddfbbo methods are both represented on the left graph and remain close to their initialization.
        Both are outperformed by the~\methodhybridfinal and~\methodhybridfinalmidterm methods.}
    \label{fig:control_results_harmonic}
\end{figure}

Figure~\ref{fig:control_results_harmonic} reports the numerical results obtained with the methods introduced in Section~\ref{sec:app_opt_control/numerical_methods} on Problem~\eqref{pb:harmonic}, with the parameters $T = 60$, $N = 600$, $M = 2$ and~$K = 1/2$. They follow our expectations from Section~\ref{sec:app_opt_control/numerical_methods}.
Among the listed methods, only the \methodhybridgeneric method performs well.
The \methodso method fails to detect the gain resulting from a large final spring elongation, and the \methoddfbbo method fails because \nomad faces too many variables.
Both methods \methodhybridfinal and \methodhybridfinalmidterm succeed because they use \nomad with few variables to detect a relevant final state and \ipopt to recover a trajectory joining it.

\subsection{A first Zermelo-type problem}
\label{sec:numerical_tests/zermelo_1}

\begin{figure}[h]
    \centering
    \includegraphics[width=\linewidth]{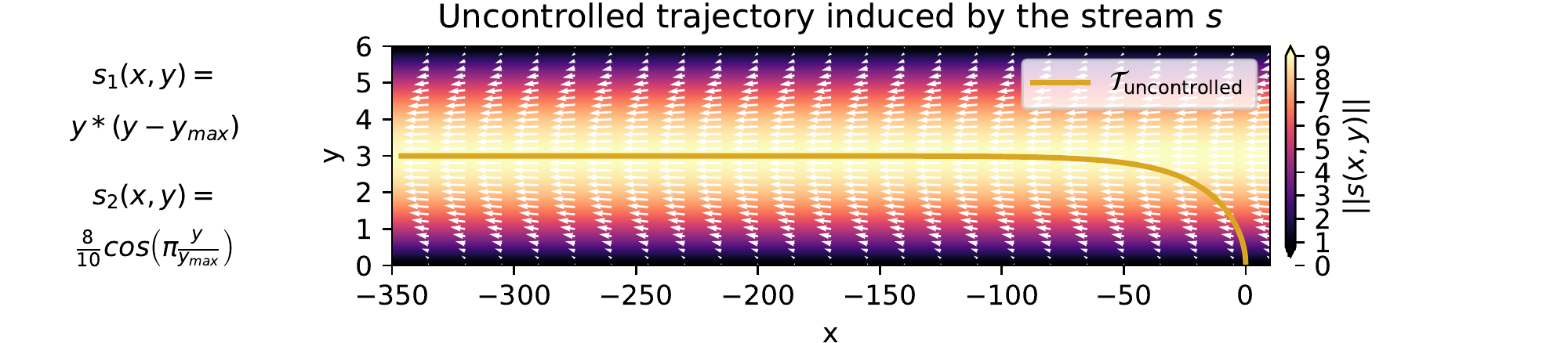}
    \caption{The stream~$s = (s_1,s_2) : \R \times[0,y_{\rm{max}}] \to \R^2$ with $y_{\rm{max}}=6$.
    It repels the boat from the shores.
        The uncontrolled trajectory (with the motor turned off, that is, with $(u,\theta)\equiv (0,0)$) goes to $x\to-\infty$, $y\to\frac{y_\mathrm{max}}{2}$.
        The horizontal component $s_1$ of the stream is powerful in the middle of the river and cannot be overcame by the driving force of the boat delivered by the motor.}
    \label{fig:zermelo_stream}
\end{figure}

In this section we consider a variant of the Zermelo problem, with a piecewise constant Mayer cost function, given by
\begin{equation}
    \tag{\mbox{$\mathcal{Z}_1$}}
    \label{pb:Zermelo_1}
    \problemoptimcontrol
        {\minimize}
        {(x,y) \in \mathcal{AC}([0,T],\R^2) \\
         (u,\theta) \in \mathcal{L}^{\infty}([0,T],\R^2)
        }
        { \displaystyle -\floor{\frac{x(T)}{10}}
         +
        \frac{1}{T} \int_0^T u(t) \, dt 
        }
        {\dot{x}(t) = s_1(x(t),y(t)) + u(t) \cos \theta(t), & \mbox{a.e.\ } t \in [0,T], \\[3pt]
        \dot{y}(t) = s_2(x(t),y(t)) + u(t) \sin \theta(t), & \mbox{a.e.\ } t \in [0,T], \\[3pt]
        x(0)=0, \\[3pt]
        y(0)=0, \\[3pt]
        0 \leq y(t) \leq y_{\rm{max}}, & \forall t \in [0,T],
        \\[3pt]
        y(T) = y_{\rm{max}}, \\[3pt]
        u(t) \in [0,1], & \mbox{a.e.\ } t \in [0,T]
        }
\end{equation}
where~$(x(t),y(t))$ is the position of a boat on a river $\R\times[0,y_{\rm{max}}]$ of width~$y_{\rm{max}} > 0$, navigating from~$(0,0)$ and trying to reach the upper shore $\R \times \{y_{\rm{max}}\}$.
The boat is equipped with a motor delivering a driving force $u(t)(\cos\theta(t),\sin\theta(t))$ with controllable norm~$u(t) \in [0,1]$ and direction~$\theta(t) \in \R$.
The river has a stream~$s = (s_1,s_2)$ (defined in Figure~\ref{fig:zermelo_stream}).
Note that the stream and the driving force actually do not represent physical forces.
Indeed the stream is a vector field of velocities and the control~$(u(t),\theta(t))$ induces an additive term to the velocity of the boat.
There are parking lots of width~$10$ in the upper shore.
The objective is to reach the rightmost parking lot at the final time, by taking into account the average value of~$u(t)$.

\begin{figure}[h]
    \centering
    \includegraphics[width=0.75\linewidth]{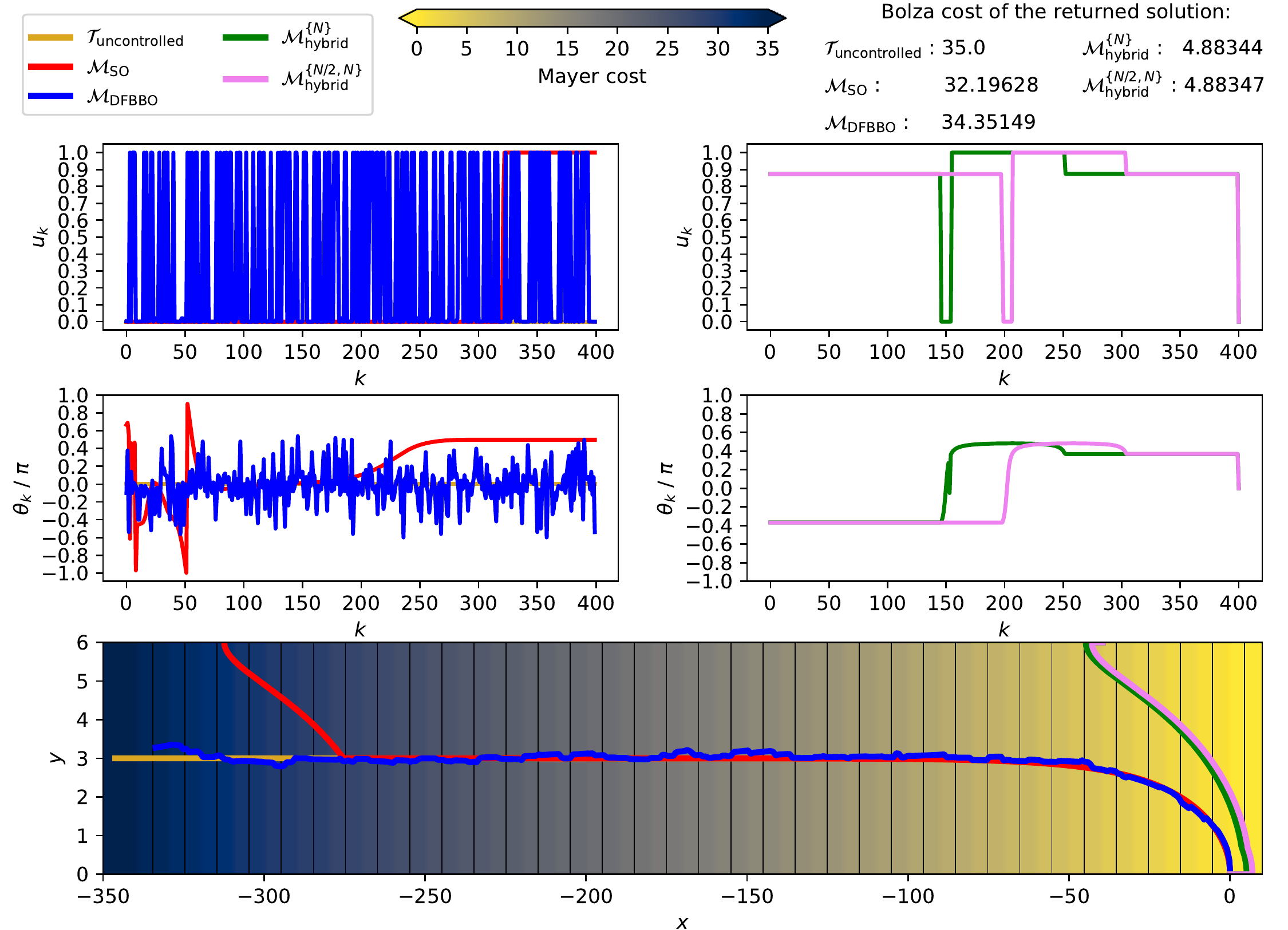}
    \caption{Numerical results obtained on Problem~\eqref{pb:Zermelo_1}.
    The \methodhybridfinal and \methodhybridfinalmidterm methods perform well.
    The \methodso method successfully reaches the feasible set but remains in a locally optimal solution, and the \methoddfbbo method fails to find a feasible trajectory.}
    \label{fig:results_zermelo_1}
\end{figure}

Figure~\ref{fig:results_zermelo_1} reports the numerical results obtained with the methods introduced in Section~\ref{sec:app_opt_control/numerical_methods} on Problem~\eqref{pb:Zermelo_1}, with the parameters $T=40$, $N=400$ and~$y_{\rm{max}}=6$.
These results follow our expectations from Section~\ref{sec:app_opt_control/numerical_methods}.
Recall that all methods start from the uncontrolled trajectory which is infeasible since~$y_N \neq y_{\rm{max}}$ (see Figure~\ref{fig:zermelo_stream}).
The \methodso first generates a feasible trajectory (with $y_N = 6$) but with a high objective value ($x_N$ is far on the negative values).
Then, it struggles to optimize it because the Mayer cost is constant on each interval $[10(q-1),10q)$ with $q \in \Z$ (up to slight variations due to the smoothing, see Remark~\ref{rq:smoothing_floor}).
This causes the~\methodso method to converge to a local optimum.
The \methoddfbbo method fails to find a feasible solution.
Both the \methodhybridfinal and \methodhybridfinalmidterm methods recover a feasible trajectory with a low objective function value.
The \methodhybridfinal and~\methodhybridfinalmidterm methods both outperform the other methods by a significant amount.

\subsection{A second Zermelo-type problem}
\label{sec:numerical_tests/zermelo_2}

In this section our objective is to highlight the expected superiority of the \methodhybridfinalmidterm method over the~\methodhybridfinal method when considering a problem with two targets~$A$ and~$B$ to approach (with respect to a discrete distance), one at the final time~$T$ and one at the midtime~$T/2$.
To this aim we consider a second variant of the Zermelo problem given by
\begin{equation}
    \tag{\mbox{$\mathcal{Z}_2$}}
    \label{pb:Zermelo_2}
    \problemoptimcontrol
        {\minimize}
        {X \in \mathcal{AC}([0,T],\R^2) \\
         (u,\theta) \in \mathcal{L}^{\infty}([0,T],\R^2)
        }
        {\floor{\norm{X(T) - A}_{\alpha}} + \floor{\norm{X(\frac{T}{2}) - B}_{\beta}} + \displaystyle \dfrac{1}{T} \int_0^T u(t)\, dt 
        }
        {\dot{x}(t) = s_1(X(t)) + u(t) \cos \theta(t), & \mbox{a.e.\ } t \in [0,T], \\[3pt]
         \dot{y}(t) = s_2(X(t)) + u(t) \sin \theta(t), & \mbox{a.e.\ } t \in [0,T], \\[3pt]
         X(0) = (0,0), \\[3pt]
         X(t) \in \R\times[0,y_{\rm{max}}], & \forall t \in [0,T], \\[3pt]
         u(t) \in [0,1], & \mbox{a.e.\ } t \in [0,T]
        }
\end{equation}
with~$X=(x,y)$, and using the control system and parameters from Section~\ref{sec:numerical_tests/zermelo_1} and where $A$, $\alpha$, $B$, $\beta \in \R^2$.

Note that the \methodhybridfinal method is theoretically unsuitable for this problem since, with notations from Section~\ref{sec:app_opt_control}, we have $\Omega = \{N/2,N\}$ and $\Lambda=\{N\}$ while the inclusion~$\Omega \subseteq \Lambda$ is required.
However, to illustrate the lack of performance of the~\methodhybridfinal method on this problem, we will implement it by adding the term $\floor{\norm{X(T/2)-B}_{\beta}}$ to the objective of Subproblem~($S(0,X_0,N,X_N)$) which is solved by \ipopt.

\begin{figure}[h]
    \centering
    \includegraphics[width=0.75\linewidth]{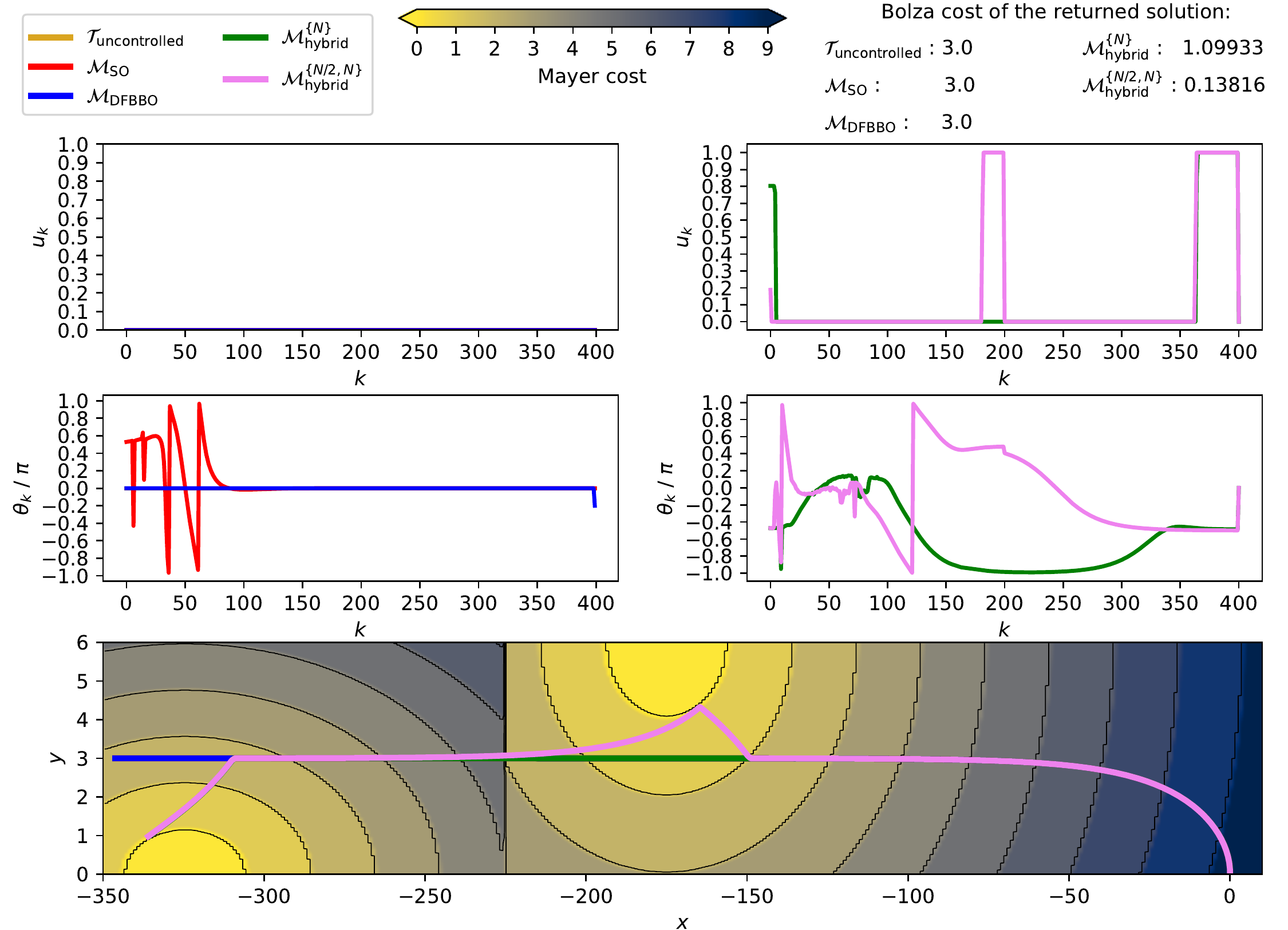}
    \caption{Numerical results obtained on Problem~\eqref{pb:Zermelo_2}.
    The \methodso and \methoddfbbo methods fail to differ from the uncontrolled trajectory.
    The \methodhybridfinal method reaches the lowest-cost plateau at the final time but fails to approach the midterm reference state.
    The \methodhybridfinalmidterm method reaches both lowest-cost plateaus.}
    \label{fig:results_zermelo_2}
\end{figure}

Figure~\ref{fig:results_zermelo_2} reports the numerical results obtained with the methods introduced in Section~\ref{sec:app_opt_control/numerical_methods} on Problem~\eqref{pb:Zermelo_2}, with the parameters $T$, $N$, $y_\mathrm{max}$ and $s$ as in Section~\ref{sec:numerical_tests/zermelo_1} and with $A = (-325,0)$, $\alpha = (\frac{1}{20},\frac{5}{6})$, $B = (-175,6)$ and $\beta = (\frac{1}{20},\frac{1}{2})$.
The plateaus associated with $B$ are represented on the lower graph on the locus~$]-225,0]\times[0,6]$, while the plateaus associated with $A$ appear on the locus~$[-350,-225{[} \times[0,6]$.

The~\methodso and~\methoddfbbo methods leave the boat uncontrolled and miss the targets.
The \methodhybridfinal and \methodhybridfinalmidterm methods perform better, although only the \methodhybridfinalmidterm method approaches the midterm target.
As expected, only the method explicitly handling both singularities is reliable.
This example shows a case where the dynamics and the Lagrange cost repel the system from a target, since the stream repels the boat from the shores while confronting it is penalized by the Lagrange cost.
Then, confronting the stream appears irrelevant for each method except the \methodhybridfinalmidterm method.

\subsection{A truncated variant of the Lotka-Volterra stabilization problem}
\label{sec:numerical_tests/lotka-volterra} 

In the three previous problems, the analysis of the numerical methods follows our expectations from Section~\ref{sec:app_opt_control/numerical_methods}.
In particular the~\methodhybridgeneric method performs better than the others, and $\Lambda$ must be chosen accordingly to the singularities appearing in the Mayer cost function.
Precisely $\Lambda$ must contain the set $\Omega$ and, at a first glance, there is no need to make it larger.
However, and surprisingly, some singular problems with $\Omega = \{N\}$ may make the \methodhybridfinalmidterm method perform better than the \methodhybridfinal method.
In this section we enlighten and analyze this phenomenon on an example, based on a simplified two-dimensional Lotka-Volterra biological system, given by
\begin{equation}
    \tag{\mbox{$\mathcal{LV}$}}
    \label{pb:lotka_volterra}
    \problemoptimcontrol
        {\minimize}
        {X = (x,y) \in \mathcal{AC}([0,T],\R^2) \\
         u \in \mathcal{L}^{\infty}([0,T],\R)
        }
        { \displaystyle \floor{\norm{X(T)-\overline{X}}_\alpha}
            +
        \frac{1}{T} \int_0^T u(t) \, dt 
        }
        {\dot{x}(t) = x(t)(y(t)-\overline{y}-u(t)), & \mbox{a.e.\ } t \in [0,T], \\[3pt]
        \dot{y}(t) = y(t)(\overline{x}-x(t)), & \mbox{a.e.\ } t \in [0,T],
        \\[3pt]
        x(0) = x_0, \\[3pt]
        y(0) = y_0, \\[3pt]
        u(t) \in [0,\frac{1}{4}], & \mbox{a.e.\ } t \in [0,T]
        }
\end{equation}
where $X(t) = (x(t),y(t))$ represents the number of thousands of individuals in a specie of predators and preys respectively,~$X_0 = (x_0,y_0) \in (\R_+^*)^2$ stands for the initial condition, $\overline{X} = (\overline{x},\overline{y}) \in (\R_+^*)^2$ denotes the equilibrium of the Lotka-Volterra system, and~$u(t) \in [0,1/4]$ is the ratio of predators decimated by an human intervention at time $t$.
With no intervention (that is, with~$u \equiv 0$), the Lotka-Volterra system has a periodical motion.
The objective of Problem~\eqref{pb:lotka_volterra} is to stabilize the system by approaching the equilibrium~$\overline{X}$ (with respect to a discrete distance), and minimizing the average human intervention.
We use the parameters $T=60$, $N=600$, $X_0 = (x_0,y_0) = (2.25,4.25)$ and $\overline{X} = (\overline{x},\overline{y}) = (1,1)$.
They lead to the final state of the uncontrolled trajectory $X_N \approx (0.05,1.59)$.

\begin{figure}[h]
    \centering
    \includegraphics[width=\linewidth]{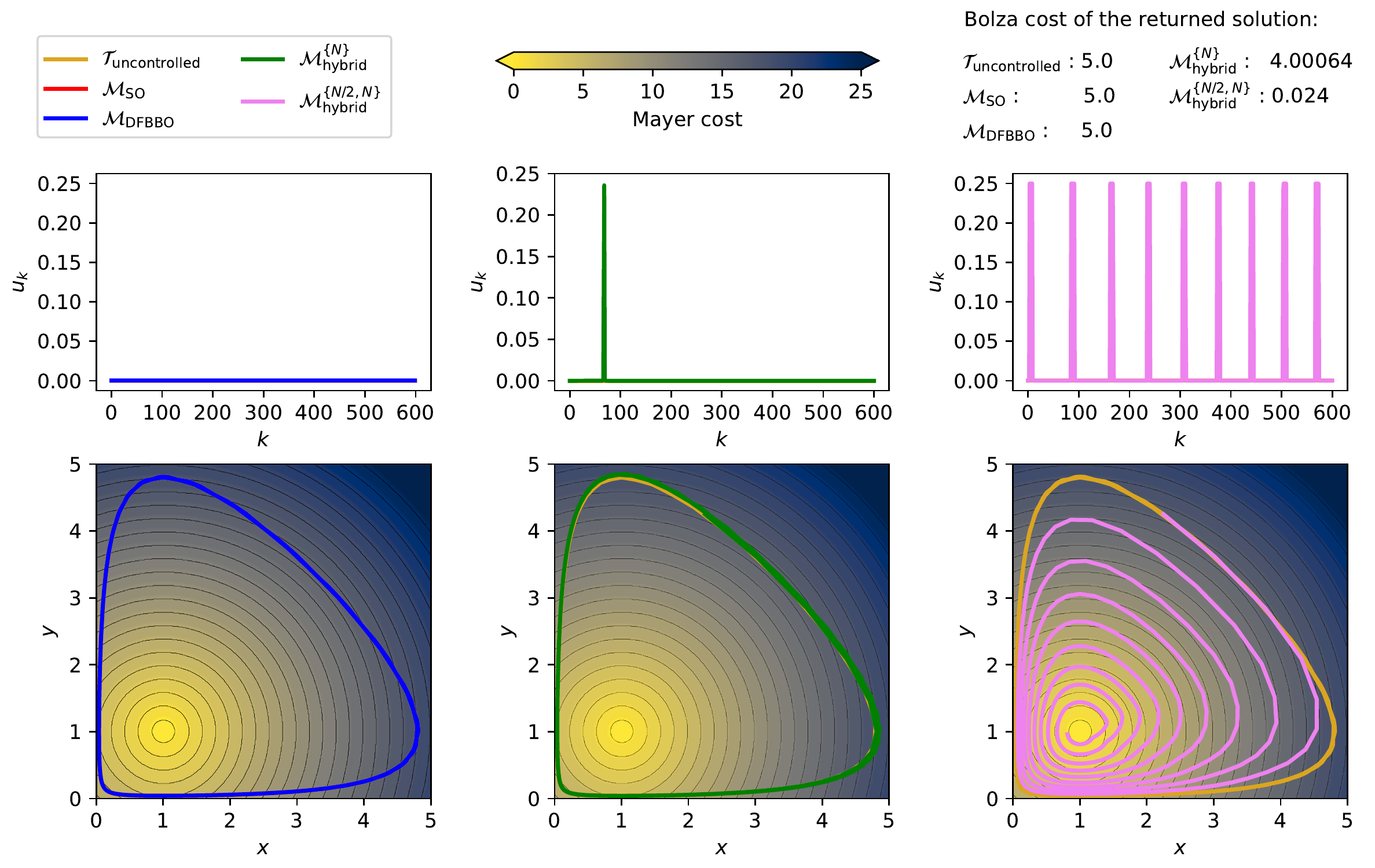}
    \caption{Numerical results obtained on Problem~\eqref{pb:lotka_volterra}.
    The \methodso and \methoddfbbo methods are both represented on the left graph and remain at their initialization.
    The \methodhybridfinal method also remains on a trajectory close to the uncontrolled one.
    The \methodhybridfinalmidterm method performs well and successfully reaches the lowest-cost plateau.}
    \label{fig:lotka_volterra_results}
\end{figure}

The controls and trajectories found by each method developed in Section~\ref{sec:app_opt_control/numerical_methods} are shown on Figure~\ref{fig:lotka_volterra_results}. The \methodso and \methoddfbbo methods fail to differ from the uncontrolled trajectory.
Contrary to the prior problems, the \methodhybridfinal method also fails.
The \methodhybridfinal method only pushes one control value at~$1/4$ while all the others remain at $0$, to slightly alter the final state~$X_N$ and make it cross one single plateau.
This behavior is intriguing.
Since this method leaves only the variable~$X_N$ to \nomad, at a first glance, there is no reason to explain why the method fails to propose a final state~$X_N$ closer to~$\overline{X}$.
This issue is of prime importance for further research works and is discussed in more details in the next paragraph. Finally note that the \methodhybridfinalmidterm method performs well since it recovers a trajectory reaching the lowest-cost plateau.

\begin{figure}[h]
    \centering
    \includegraphics[width=0.60\linewidth]{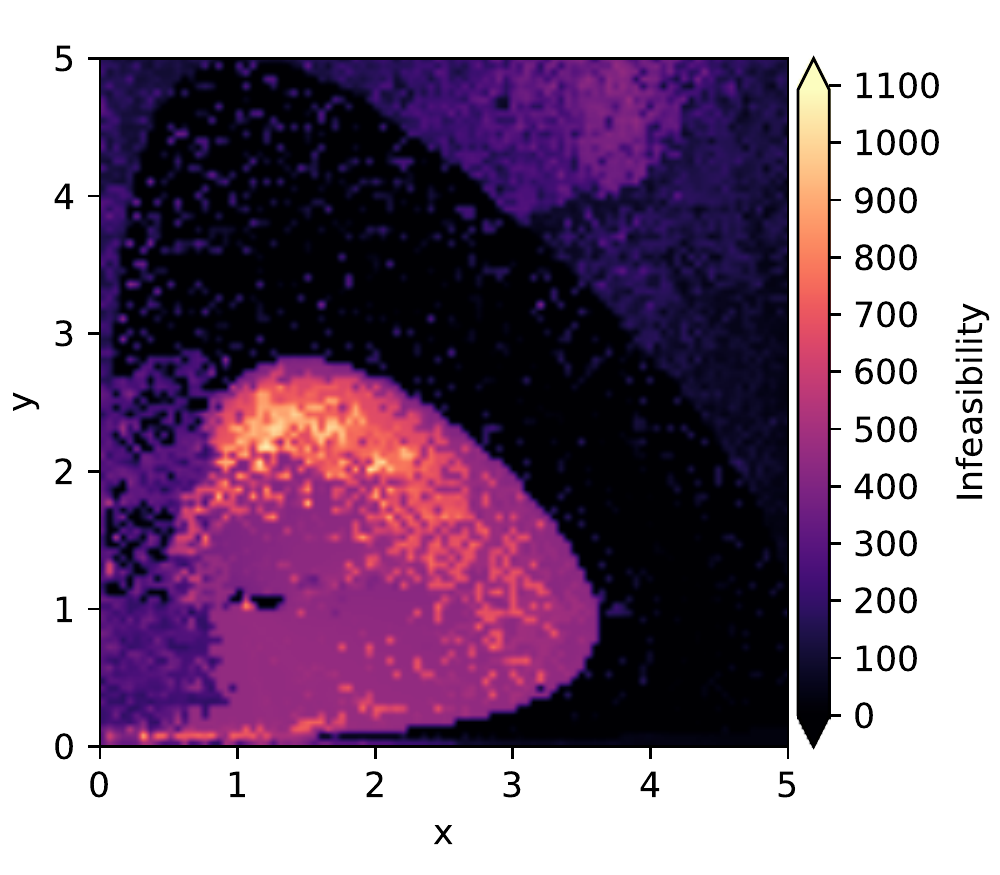}
    \caption{The infeasibility of the solution returned by \ipopt attempting to solve Subproblem~$(\mathrm{S}(0,X_0,N,X_N))$ for each~$X_N = (x,y) \in [0,5]^2$.
        It may differ from the true minimal infeasibility because \ipopt may converge to a locally infeasible solution in its attempt to solve a feasible Subproblem~$(\mathrm{S}(0,X_0,N,X_N))$.}
    \label{fig:lotka-volterra_infeasibility}
\end{figure}

As suggested by the numerical results obtained, Problem~\eqref{eq:Rcontrol_xm} is difficult to solve numerically in the context of this section.
Recall that, in Problem~\eqref{eq:Rcontrol_xm}, the only variable is the final state $X_N$, while the search of a trajectory linking $X_0$ to~$X_N$ is performed via Subproblem~$(\mathrm{S}(0,X_0,N,X_N))$.
Recall also that \ipopt looks for a (feasible and) optimal solution of the subproblem, but it has no guarantee to localize neither an optimal solution nor a feasible one.
As it aims primarily a feasible solution before optimizing it, \ipopt may converge to an infeasible solution in its attempt to solve a feasible subproblem, because this solution is a local minimum of the infeasibility measure.
Hence the variable $X_N$ may trigger an hidden constraint in Problem~\eqref{eq:Rcontrol_xm}, as it may generate a subproblem for which \ipopt fails to find a feasible solution.
This constraint does not appear on the mathematical formulation of Problem~\eqref{pb:lotka_volterra}, as it comes from the computer code implementing the \methodhybridgeneric method, but it has an impact on the returned solution.
This phenomenon explains why \nomad struggles to solve Problem~\eqref{eq:Rcontrol_xm} in the context of this section.
Indeed, consider Figure~\ref{fig:lotka-volterra_infeasibility} which shows the infeasibility of the solution returned by \ipopt on Subproblem~($\mathrm{S}(0,X_0,N,X_N)$) for each point $X_N$.
The set of $X_N$ leading to a subproblem for which the solution returned by \ipopt has a high infeasibility measure surrounds the equilibrium~$\overline{X}$.
Therefore, starting from the final state of the uncontrolled trajectory~$X_N \approx (0.05,1.59)$ in the large black ribbon of Figure~\ref{fig:lotka-volterra_infeasibility}, \nomad sees its initialization as a feasible solution but also sees the colored zone as a zone of high infeasibility, and thus remains in the black ribbon.
Finally, as the \methodhybridfinalmidterm method performs better than the \methodhybridfinal method, we suspect that Subproblem~$(\mathrm{S}(0,X_0,N,X_N))$ recovering the entire trajectory is harder to solve numerically by \ipopt than the two Subproblems~$(\mathrm{S}(0,X_0,N/2,X_{N/2}))$ and~$(\mathrm{S}(N/2,X_{N/2},N,X_N))$ recovering half of the trajectory each.
Furthermore we even found that \ipopt fails to find a feasible solution of the feasible Subproblem~$(\mathrm{S}(0,X_0,N,X_N))$ taking as~$X_N$ the final state of the trajectory returned by the~\methodhybridfinalmidterm method.
This subproblem is feasible (since the trajectory returned by the~\methodhybridfinalmidterm method is a feasible solution), but \ipopt fails to detect a feasible solution.

\section{Discussion}
\label{sec:conclusion}

This work introduces an hybrid numerical method, relying on two algorithms from \dfo and \so, to solve optimal control problems with a piecewise constant Mayer cost function.
This hybrid method is based on a reformulation of the discretized problem which splits the search of a full trajectory-control pair into several subproblems searching for partial trajectory-control pairs.
The state variables that may impact the singularities of the Mayer cost function are optimized on the main problem via a \dfo algorithm, while the search of a trajectory joining those is left to smooth subproblems solved via a \so algorithm.

As evidenced in Section~\ref{sec:numerical_tests}, the \methodhybridgeneric method solves efficiently optimal control problems with a piecewise constant Mayer cost function, while non-hybrid approaches fail.
Although we voluntarily restricted our work to optimal control problems with a piecewise constant Mayer cost function only, our method can be extended to a larger class of optimal control problems, including for example nondifferentiable cost functions, piecewise continuous cost functions, extended-real valued cost functions, etc., or cases where the control variables may also impact the singularities.
From a general point of view, any cost or constraint function, or any dynamics, involving a singularity handled by \dfo algorithms may be handled by an adapted hybrid method.
However, this last class of problems is large and some other methods may be more appropriate for some of its subclasses.
A future development could rely in a better identification of the class of problems on which hybrid methods work well.

As explained above, the \methodhybridgeneric method can be extended in several directions. It can also be improved in several aspects.
First, we know from Section~\ref{sec:numerical_tests/lotka-volterra} that a relevant choice of the subset~$\Lambda$ does not simply consist in making it as small as possible, even if a subset~$\Lambda$ with small cardinality is preferred.
Hence, a possible improvement is to fix the cardinal of $\Lambda$ and to leave its elements as variables of the main problem. Note that such a process would not require any algorithm from combinatorial optimization, since \dfo algorithms handle problems with both discrete and continuous variables~\cite{AuLeDTr2018}.
Preliminary numerical experiments are promising and this extension is expected to be the topic of future research works.

Second, we noticed that the feasibility of the subproblems highly depends on the values of the states optimized on the main problem, to the point that it may be difficult to find states leading to feasible subproblems with a \dfo algorithm (see Figure~\ref{fig:lotka-volterra_infeasibility} for an example).
A possible improvement would be to relax the constraint forcing the subproblems to match exactly the states proposed by the main problem.
The trajectory-finding subproblems may be relaxed to simply approach the states instead of matching them.
This could be done by adding a penalization term to the objective functions of the subproblems, penalizing the distance between the states of the system at the times indexed by $\Lambda$ and the proposed states at these times.

Third, the known problem structure may be exploited to guide the main \dfo algorithm.
For example, because of the dynamics of the control system constraining the trajectory, forcing the system to reach a given state at a given time restricts the space of admissible states at future times.
Hence, in the~\methodhybridgeneric method, there are links between the state variables optimized at the main level.
Furthermore, we may construct surrogates~\cite{BoDeFrSeToTr99a,ConnScheinbergVicente,CustodioScheinbergVicente,sgtelib} of the subproblems to obtain a fast approximation of their optimal solutions.
For example, we may launch a \so algorithm and stop its execution after a few iterations, and consider the low-precision solution.
The \dfo algorithm may use this procedure to estimate how relevant a potential solution is, with no need to systematically solve the subproblems with full accuracy.
Such a development could significantly reduce the computational time of the \methodhybridgeneric method.

\bibliographystyle{unsrtnat}


\end{document}